\documentclass[11pt]{amsart}
\usepackage[utf8]{inputenc}
\usepackage[T1]{fontenc}
\usepackage{graphicx,amssymb,amsfonts,latexsym,amsmath,amsthm, times, amscd, euscript, MnSymbol}
\usepackage{mathrsfs}
\usepackage{epsfig}
\usepackage{pst-text}
\usepackage{enumerate}
\usepackage{color}
\usepackage[english]{babel}

\numberwithin{equation}{section}
\usepackage[active]{srcltx}

\usepackage{verbatim}
\usepackage{multicol}
\usepackage{bezier}
\usepackage[retainorgcmds]{IEEEtrantools}

\usepackage[running]{lineno}    

\newtheorem{theorem}{Theorem}[section]
\newtheorem{lemma}[theorem]{Lemma}
\newtheorem{corollary}[theorem]{Corollary}
\newtheorem{proposition}[theorem]{Proposition}
\newtheorem{definition}[theorem]{Definition}
\newtheorem{conjecture}[theorem]{Conjecture}

\newtheorem{problem}{\noindent {\bf Problem}}

\newtheorem{claim}{\noindent {\bf Claim}}
\newtheorem{question}{\noindent {\bf Question}}

\newcommand{\NN}{{\mathbb N}}

\newcommand{\age}{\mathcal{A}}

\def\endproof{\hfill {\kern 6pt\penalty 500
\raise -0pt\hbox{\vrule \vbox to5pt {\hrule width 5pt
\vfill\hrule}\vrule}}}

\title[Monomorphic decomposition and profile]{Ordered structures with no finite monomorphic decomposition. Application to the profile of hereditary classes.}

\author[D.Oudrar] {Djamila Oudrar*}
\thanks{*The author was supported by CMEP-Tassili grant}
\address{Faculty of Mathematics, USTHB, Algiers, Algeria}
\email {dabchiche@usthb.dz; djoudrar@gmail.com}
\author [M.Pouzet]{Maurice Pouzet}
\address{Univ. Lyon, Universit\'e Claude-Bernard Lyon1, CNRS UMR 5208, Institut Camille Jordan,  43 bd. 11 Novembre 1918, 69622 Villeurbanne Cedex, France and University of Calgary, Department of Mathematics and Statistics, Calgary, Alberta, Canada T2N 1N4} \email{pouzet@univ-lyon1.fr}

\date{\today}


\begin{document}

\maketitle

\dedicatory{\dagger Dedicated to  the memory of Roland Fra\"{i}ss\'e and Claude Frasnay}

\begin{abstract}
We present a structural approach of some results about jumps  in the  behavior of the profile (alias generating function) of hereditary classes of finite
structures. We consider  the following notion due to N.Thi\'ery and the second author.
A \emph{monomorphic decomposition} of a relational structure $R$ is a  partition of its domain $V(R)$ into a family of sets $(V_x)_{x\in X}$ such that the
restrictions of $R$ to two finite subsets $A$ and $A'$  of $V(R)$ are isomorphic provided that the traces $A\cap V_x$ and $A'\cap V_x$ have the same size for each
$x\in X$. Let  $\mathscr S_\mu $ be the class of relational structures of signature $\mu$ which do not have a finite monomorphic decomposition. We show that if a
hereditary subclass $\mathscr D$ of $\mathscr S_\mu $ is made of  ordered relational structures then it contains a finite subset $\mathfrak A$ such that every
member of $\mathscr D$ embeds  some member of $\mathfrak A$. Furthermore,  for each $R\in \mathfrak A$ the  profile of  the  age $\age(R)$  of $R$ (made of
finite substructures of $R$) is at least exponential.   We  deduce  that if the profile of a hereditary class of finite ordered structures  is not bounded above by a
polynomial then it is at least exponential. For ordered graphs, this result is a part of  classification obtained by  Balogh, Bollob\'as and Morris (2006).
\end{abstract}

\section{Introduction and presentation of the results}
The   \textit{profile} of a class  $\mathscr {C}$  of finite relational structures is the integer function  $\varphi _{\mathscr{C}}$  which counts for each non negative integer $n$  the number of members of   $\mathscr{C}$
on $n$ elements, isomorphic structures being  identified.  The behavior of this function has been discussed in many papers, particularly
when  $\mathscr{C}$ is
\emph{hereditary} (that is contains  every substructure of any member of  $\mathscr{C}$) and is made of graphs (directed or not),  tournaments, ordered sets,  ordered  graphs or ordered hypergraphs. Futhermore, thanks to a result of Cameron \cite{cameron},  it turns out that the line of study about permutations (see \cite{A-A}) originating in the Stanley-Wilf conjecture, solved by   Marcus and Tard\"os (2004) \cite{Marcus},
falls under the frame of the profile of hereditary classes of ordered relational structures (see \cite{oudrar-pouzet1, oudrar-pouzet2}).
The results show that the profile cannot be arbitrary: there are jumps in its possible growth rate. Typically, its growth  is polynomial or faster than every polynomial (\cite {pouzet.tr.1978} for ages, see  \cite{pouzet-profile2006} for a survey) and for several classes of structures, it is either at least exponential (e.g. for tournaments \cite{BBM07,Bou-Pouz}, ordered graphs and hypergraphs \cite{B-B-M06,B-B-M06/2,klazar1} and permutations \cite{kaiser-klazar}) or  at least the growth of the partition function (e.g. for graphs \cite {B-B-S-S}). For more, see  the survey of Klazar \cite{klazar}.

In this paper,  we consider hereditary classes of ordered relational structures.  We describe  those with polynomially bounded profile; we prove that  each hereditary class  whose profil is not bounded by a polynomial contains a minimal one (minimal with respect to  inclusion) and prove that its profile is exponential. The case of ordered binary relational structures and particularly the case of ordered  irreflexive directed graphs was treated in Chapter 8 of \cite{oudrarthese} and are presented in \cite{oudrar}.  On the surface, the  case of general  ordered binary relational structures is similar, but the case  of ternary relations is more involved.

Let us present our main results. Each structure we consider is of the form $R:= ( V, \leq, (\rho_{j})_{j\in J})$, where $\leq$ is a linear order on $V$ and each  $\rho_j$ is a \emph{$m_j$-ary relational structure on $V$}, that is a subset of $V^{m_j}$ for some non-negative integer $m_j$, the \emph{arity} of $\rho_j$. We will say that the  sequence $\mu:= (m_j)_{j\in J}$ is the  \emph {restricted signature} of $R$. The \emph {age} of $R$ is the set $\age(R)$ consisting of the structures induced by $R$ on the  finite subsets of $V$, these structures being considered up to isomorphy. A relational structure of the form $(V, \leq)$ where $\leq$ is a  linear order on $V$ is  a \emph{chain};  if it is  of the form $B:= (V, \leq, \leq')$ where $\leq$ and $\leq'$ are two linear orders on $V$ this is a \emph{bichain}, and if it is of the form $G: =(V, \leq, \rho)$ where $\rho$ is a binary relation this is an \emph{ordered directed graph}. Chains, bichains and ordered directed graphs are the basic examples of ordered structures.

An \emph{interval decomposition} of $R$ is a partition $\mathcal P$ of $V$ into  intervals $I$ of the chain $C:= (V, \leq)$ such that for every integer $n$ and every pair $A,~A'$ of $n$-element subsets of $V$, the induced structures on $A$ and $A'$ are isomorphic whenever the traces $A\cap I$ and  $A'\cap I$ have the same number of elements for each interval $I$.  For  example, if $R$ is  the bichain $(V, \leq, \leq')$, $\mathcal P$ is an interval  decomposition of $V$ iff each block $I$ is an interval for each of the two orders and they coincide or are opposite on $I$ (see \cite{monteil-pouzet}).

If  a relational structure $R$ has  an interval decomposition into  finitely many blocks, say $k+1$, then trivially, $\varphi _{\age(R)}$, the profile of $\age(R)$, is bounded above by some polynomial whose degree is, at most, $k$. According to Theorem 1.7 of  \cite{P-T-2013},  this is eventually a \emph{quasi-polynomial}, that is for $n$ large enough, $\varphi _{\age(R)}(n)$ is a sum $a_{k}(n)n^{k}+\cdots+ a_0(n)$ whose coefficients $a_{k}(n), \dots, a_0(n)$ are periodic functions. In this paper,   we  show  that the profil  is eventually  a polynomial (Theorem \ref{thm: polynomial-interval}).

We prove that, essentially,  the converse holds.

 \begin{theorem}\label{thm:poly-expo1} Let  $\mathscr C$ be a hereditary class of finite ordered relational structures with a finite restricted signature $\mu$. Then,   either there is  some integer $k$ such that every member of $\mathscr  C$ has an interval decomposition into at most $k+1$ blocks, in which case  $\mathscr C$ is a finite union of ages of ordered relational structures, each having an interval decomposition into at most $k+1$ blocks,  and the profile of $\mathscr C$ is eventually a polynomial (of degree at most $k$), or the profile of $\mathscr C$ is at least  exponential.
\end{theorem}
\medskip

The jump of the growth of profile from polynomial to exponential  was obtained for bichains by Kaiser and Klazar \cite{kaiser-klazar} and extended  to ordered graphs by  Balogh,  Bollob\'as and Morris (2006) (Theorem 1.1  of \cite{B-B-M06/2}). Their  results go much beyond exponential profile.

\medskip

The  proof of Theorem \ref{thm:poly-expo1} has three steps that we indicate now.

The first  is a reduction to the case where $\mathscr C$ is of the form $\age(R)$. This follows from    the next lemma.

\begin{lemma} \label{lem:reduction}
If a hereditary class $\mathscr C$  of finite ordered relational structures with a finite signature  contains for every integer $k$ some finite  structure which has no interval decomposition into at most $k+1$ blocks, then it contains a hereditary class  $\mathscr A$ with the same property which is minimal w.r.t. inclusion. This class is the age  of some ordered relational structure $R$ which does not have any finite interval decomposition.

\end{lemma}

The proof of Lemma  \ref{lem:reduction} relies on properties of well-quasi-ordering and of ordered structures.
A more general result is given in Lemma \ref{lem:reductionbis}.   The concepts behind the proof of Lemma  \ref{lem:reduction}  are important enough to be presented now.
To prove  that ${\mathscr C}$ contains  a  minimal hereditary subclass  ${\mathscr A}$ that contains for every integer $k$ some finite  structure which has no interval decomposition into at most $k+1$ blocks,  we prove that ${\mathscr C}$  contains a hereditary subclass ${\mathscr C'}$ with this property and with no  infinite antichain (w.r.t. embeddability). If  ${\mathscr C}$ contains some  infinite  antichain, then  a straightforward application of Zorn's Lemma yields that among the final segments  of ${\mathscr C}$ which have infinitely many minimal elements,  there is some, say $\mathcal F$, which is maximal w.r.t. inclusion. By construction, $\mathscr C':={\mathscr C}\setminus \mathcal F$ has no infinite antichain (for a detailed argument, see  Theorem \ref{minimalage}). Furthermore, $\mathscr C'$ contains for every integer $k$ some finite  structure which has no interval decomposition into at most $k+1$ blocks.
The proof of this  fact is deeper. Indeed, if does not hold,   an old  result of the second author about hereditary well-quasi-ordered classes implies that  $\mathscr C'$ would have  finitely many bounds  (indeed,  $\mathscr C'$ would be hereditary w.q.o.  and hence would have finitely many bounds, see Proposition \ref{prop:hw.q.o.}), which is not the case. Now, if a hereditary class contains no infinite antichain, then by a result of Higman \cite{higman.1952} about well quasi orders, the  collection of hereditary subclasses is well founded w.r.t. inclusion, hence  we may find a hereditary class  $\mathscr A$ with the same property as ${\mathscr C}$ which is minimal w.r.t. inclusion. Clearly, $\mathscr A$ cannot be the union of two proper hereditary classes, hence it must be  up-directed w.r.t. embeddability. Thus,  according to an old and well know result of Fra\"{\i}ss\'e \cite{fraisse} p.279,  this is the age of a relational structure $R$.  This structure is ordered and does not have any finite interval decomposition.

Hence our reduction is done.

For the second step,
 we  introduce the class $\mathscr O_{\mu}$  of ordered  relational structures of restricted signature $\mu$, $\mu$ finite, which do not have a finite interval decomposition. We prove:


 \begin{theorem}\label{theo:base}
There is a finite subset $\mathfrak A$ made of incomparable structures of $\mathscr O_{\mu}$ such that every member of $\mathscr O_{\mu}$ embeds some member of $\mathfrak A$.
\end{theorem}

For example,  if  $\mathscr D$ is the subclass of  $\mathscr O_{\mu}$ made of bichains then  $\mathfrak A\cap \mathscr D$ has twenty elements \cite{monteil-pouzet}. If $\mathscr D$ is made of ordered reflexive (or irreflexive) directed graphs, $\mathfrak A\cap \mathscr D$ contains  $1242$ elements (cf Theorem 1.6 of \cite {oudrar}).

The proof, given in Section \ref{section:proofs} relies on Ramsey's theorem and the notion of monomorphic decomposition.  This  notion
was   introduced in \cite{P-T-2013} in the sequel of  R.~Fra\"{\i}ss\'e who invented the notion of monomorphy  and C.~Frasnay who proved the central result about this notion \cite{frasnay 65}.
We define an equivalence relation $\equiv_R$ on the domain of a relational structure $R$ (ordered or not) whose classes form a monomorphic decomposition.  This equivalence relation is an intersection of equivalences relations $\equiv_{k, R}$. When $R$ is ordered, we note  that every interval decomposition  of $R$ is a monomorphic decomposition and that each infinite equivalence class of  $\equiv_R$ is an interval of $R$ (Theorem \ref{thm:intervalsecomposition}). Furthermore, there  is some integer $k$ such that
 $\equiv_R$ and $\equiv_{k,R}$ coincide (Theorem \ref{thm:marytreshold}).

Using that notion, we prove in Theorem \ref{mainthm} that the members of $\mathfrak A$ have a special form: there are almost-multichainable (a notion introduced by the second author in his th\`ese d'\'Etat \cite{pouzet.tr.1978} which appeared in \cite{pouzet-impartible} and \cite{pouzet-minimale}).

In the third step,  we prove in Section \ref{section:proofs} that   the profile of members of  $\mathfrak A$  grows at least exponentially.
Let  $k\in N$ and $\mathscr O_{\mu, k}$ be  the class of ordered structures  $R$ of restricted signature $\mu$ such that the equivalence relation $\equiv_{\leq k,R}$ has infinitely many classes.

\begin{theorem}\label{theo:expo} If $R\in \mathscr O_{\mu,k}$  then the profile $\varphi_R$ is at least exponential.  Indeed, for $n$ large enough it  satisfies $\varphi_R(n) \geq d.c^{n}$ where $c$ is the largest solution of $X^{k+1}-X^{k}-1=0$ and $d$ is a positive constant depending upon $k$.
\end{theorem}
In the case of ordered undirected graphs, it was proved in \cite {B-B-M06/2}  that this profile  grows as fast as the Fibonacci sequence.

 From Lemma \ref{lem:reduction} and Theorem  \ref{thm:poly-expo1}, we can deduce the following.

\begin{corollary} \label{cor:minimal}If the profile of a hereditary class $\mathscr C$ of finite ordered relational structures (with a finite restricted signature $\mu$) is not bounded above by a polynomial then it contains a hereditary class  $\mathscr A$  which is minimal w.r.t. inclusion among the classes having this property.
\end{corollary}

\begin{proof}
If   the profile of $\mathscr C$ is not bounded above by a polynomial then  for every integer $k$, ${\mathscr C}$ contains some finite  structure which has no interval decomposition into at most $k+1$ blocks. According to Lemma \ref{lem:reduction} it contains a hereditary class  $\mathscr A$ with this property which is minimal w.r.t. inclusion. As already observed, $\mathscr A$ is the age of some relational structure. According to Theorem \ref{thm:poly-expo1}, the profile of $\mathscr A$ is exponential. If $\mathscr A'$ is a proper subclass of $\mathscr A$, its members have  an interval decomposition of bounded length, hence as observed, the profile of $\mathscr A'$ is  bounded by a polynomial. Thus $\mathscr A$ is minimal among the classes whose profil is not bounded by a polynomial. The proof of Corollary \ref{cor:minimal} is complete.
\end{proof}

This result holds for arbitrary hereditary classes of relational structures, provided that their arity is finite, (Theorem \ref{thm:minimal1}).   It appears in a somewhat equivalent form as Theorem 0.1 of \cite{P-T-2013}. It is not trivial (the argument used  in Corollary \ref{cor:minimal} does apply). The key  argument relies on a result going back to the thesis of the second author \cite{pouzet.tr.1978}, namely Lemma 4.1 p. 23 of  \cite{P-T-2013}. No complete proof has been yet published. A hint is given in Subsection \ref{subsection:profile}.


Our results  are part  of the study  of monomorphic decompositions of  a  relational structure.
The profile of a class of finite relational structures, not necessarily ordered, each  admitting  a finite monomorphic decomposition in at most $k+1$ blocks,  is the union of finitely many ages of relational structures admitting a finite monomorphic decomposition in at most $k+1$ blocks  and is  eventually a quasi-polynomial (this is the main result of \cite{P-T-2013}). Lemma \ref{lem:reduction} extends, but we do not know if Theorem \ref{theo:base} extends in general. We state that as a conjecture.

Let $\mathscr S_{\mu}$ be the class of all relational structures of signature $\mu$, $\mu$ finite,  without any finite monomorphic decomposition.

\begin{conjecture}
The class $\mathscr S_{\mu}$  contains  a finite subset $\mathfrak A$ made of incomparable structures such that every member of $\mathscr S_{\mu}$ embeds some member of $\mathfrak A$.
\end{conjecture}

As this will become apparent  in Section \ref{subsection: hypomorphy},  the difficulty is  with ternary structures.
We may note that if one  restricts  $\mathscr S_{\mu}$  to tournaments, there is a set $\mathfrak A$ with  twelve elements \cite{Bou-Pouz}. The first author has shown that the conjecture holds if $\mathscr S_{\mu}$ consists of ordered binary structures. She proved that for ordered  directed reflexive  graphs,  $\mathfrak A$  contains $1242$ elements \cite{oudrar}.  We show that if we consider the  class of undirected graphs, $\mathfrak A$ has ten elements. The proof is easy, we give it in order  to illustrate in a simple setting the technique used in the proof of Theorem \ref{theo:base}, even that graphs are not ordered structures.   We may note that a graph has a finite monomorphic decomposition iff it decomposes into a finite lexicographic sum of cliques and independent sets. Hence our later result can be stated as follows:

\begin{proposition}\label{prop: graphs}There are ten infinite undirected graphs such that an undirected  graph does not decompose into a finite lexicographical sum of cliques and independent sets iff it contains a copy of one of these ten graphs.
\end{proposition}

We may note that some of these graphs have polynomial profile, hence our machinery is not sufficient to illustrate the jump in profile beyond polynomials.

Results of this paper are included in Chapter 7 of the thesis of the first author \cite {oudrarthese}. Some have been included in a paper  posted on ArXiv \cite{oudrar-pouzet2014} and presented at ICGT 2014 (June 30-July 4  2014, Grenoble) and also at ALGOS 2020 (August 26-28, 2020, Nancy). We are pleased to thanks the organizers for these beautiful meetings and for giving us the opportunity of presenting our work.

\section{Interval decomposition and monomorphic decomposition}
In this section, we  prove an extension of Lemma \ref{lem:reduction} and Corollary \ref{cor:minimal} (Lemma \ref{lem:reductionbis} and Theorem \ref{thm:minimal1}) and the detail of the first part of Theorem  \ref{thm:poly-expo1} (Theorem \ref{thm: polynomial-interval}).  For this purpose, we recall the basic notions of the theory of relations (subsection \ref{subsection:basic}) and the notion of monomorphic decomposition with the main results of \cite{P-T-2013} (subsection \ref{subsection:monomorphic}). One of our tools is the notion of well-quasi-order (subsection \ref{subsection:w.q.o.}).

\subsection{Embedabbility, ages and  hereditary classes}\label{subsection:basic}
Our terminology follows Fra\"{\i}ss\'{e} \cite{fraisse}. Let $m$ be a non-negative integer. A \emph{$m$-ary relation} with \emph{domain} $V$ is a subset $\rho$ of $V^m$; when needed, we identify $\rho$ with its characteristic function. A \textit{relational structure}  is a pair $R:= (V,(\rho_i)_{i\in I})$ where each $\rho_i$ is a $m_i$-ary relation on $V$; the set $V$ is the \emph{domain} of $R$, denoted by $V(R)$; the sequence $\mu:=(m_i)_{i\in I}$  is the \emph{signature} of $R$.  This is a  \emph{binary relational structure}, \emph{binary structure} for short, if it is made only of binary relations. It is \emph{ordered} if we may write $I=J\cup \{0\}$ with $0\not \in J$ and if $\rho_0$ is a linear order on $V$. In simpler terms, $R= (V, \leq, (\rho_j)_{j\in J})$ where $\leq$ is a linear order on $V$. The sequence $(m_j)_{j\in J}$ is the \emph{restricted signature} of $R$ (this is  the signature of $(V, (\rho_j)_{j\in J})$).
We say that the signature is \emph{finite} if its domain $I$ is finite.
Let $R:= (V,(\rho_i)_{i\in I})$ be  a relational structure; the \emph{substructure induced by $R$ on a
  subset $A$ of $V$}, simply called the \emph{restriction of $R$ to
  $A$}, is the relational structure $R_{\restriction A}:= (A,
(A^{m_i}\cap \rho_i)_{i\in I})$. The restriction to $V\setminus \{x\}$ is denoted $R_{-x}$.  The notion  of \emph{isomorphism} between relational structures
is defined in the natural way. A
\emph{local isomorphism} of $R$ is any isomorphism between two restrictions
of $R$.  A relational structure $R$ is \emph{embeddable} into a relational structure $R'$,  in notation $R\leq R'$,  if $R$ is isomorphic to some restriction of $R'$. Embeddability   is a quasi-order on the class of structures having a given arity. We denote by $\Omega_\mu$ the class  of \emph{finite} relational structures of signature $\mu$,   quasi-ordered by embeddability. Most of the time we consider its members up to isomorphy. The \emph{age} of a relational  structure $R$ is the
set $\age(R)$ of restrictions of $R$ to finite subsets of its domain, these restrictions being considered up to isomorphy.  A class $\mathscr{C}$ of structures is \textit{hereditary} if it contains every relational structure which can be embedded in some member of $\mathscr{C}$.  If $\mathscr{B}$ is a subset of $\Omega_\mu$ then $Forb(\mathscr{B})$ is the subclass of members of $\Omega_\mu$ which embed no member of $\mathscr{B}$. Clearly, $Forb(\mathscr{B})$ is an hereditary class. For the converse, note that  a \textit{bound} of a hereditary class $\mathscr{C}$ of finite relational structures is every $R$ not in $\mathscr{C}$ such that every $R'$ which strictly embeds into $R$ belongs to $\mathscr{C}$. Clearly, bounds are finite and  if $\mathscr{B(C)}$ denotes the collection of  bounds of $\mathscr{C}$ then $\mathscr{C}=Forb(\mathscr{B(C)})$. The above notions can be expressed in purely order theoretic terms: a hereditary class is any initial segment of the (quasi) ordered set $\Omega_\mu$; each age is an ideal, that is a non-empty  initial segment which is up-directed (the converse holds if $\mu$ is finite, (Fra\"{\i}ss\'e, 1954, see \cite{fraisse}, p.279)); the bounds of a hereditary class $\mathscr C$ are the minimal elements of $\Omega_\mu\setminus \mathscr C$.


\subsection{Reducts, chainability and monomorphy}
 Let $R:= (V,(\rho_i)_{i\in I})$ be  a relational structure; if $I'$ is a subset of $I$ then $R^{I'}:= (V,(\rho_i)_{i\in I'})$ is a \emph{reduct} of $R$;  \emph{a finite reduct} if $I'$ is finite. Let $S$ be an other relational structure with the same domain as $R$. We say that $R$ is \emph{freely interpretable}  by $S$ if every local isomorphism of $S$ is a local isomorphism of $R$.  If $S$ is a chain, we say that $S$ \emph{chains} $R$.
We say that $R$ is \emph{chainable} if some chain $C$ chains $R$.
 Let $p$ be a non-negative integer. A relational structure $R$ is \emph{$p$-monomorphic}  if its restrictions to finite sets of the same cardinality $p$ are isomorphic. The structure $R$ is \emph{monomorphic} if it is $p$-monomorphic for every integer $p$.  Since two finite chains with the same cardinality are isomorphic, chains are monomorphic and hence chainable relational structures too. Conversely,  every infinite monomorphic relational structure is chainable (Fra\"{\i}ss\'e, 1954). In his thesis, Frasnay 1965 (\cite{frasnay 65})
obtained the existence of an integer $p$ such that every $p$-monomorphic relational structure $R$ whose maximum of the signature is at most $m$ and domain infinite or sufficiently large is chainable.

 The least integer $p$ in the sentence above is the  \emph{monomorphy treshold}, $d_m$.
 Its  value was given by Frasnay in 1990  \cite{frasnay 90}.
 \begin{theorem} $d_1=1$, $d_2= 3$, $d_m=2m-2$ for $m\geq 3$.
 \end{theorem}
 For a detailed exposition of this result,  see \cite {fraisse} Chapter 13, notably  p. 378.

 Chainability and monomorphy were generalized as follows. Let $F$ be a subset of $V$. We say that $R$ is \emph{$F$-chainable}  if there is a chain $C$ on $V\setminus F$ such that every local isomorphism of $C$ extended by the identity on $F$ is a local isomorphism of $R$; we say that  $R$  is \emph{$F$-monomorphic} if for every integer $n$ and every pair of $n$-element subsets $A$ and $A'$ of $V\setminus F$, there is an isomorphism of  $R_{\restriction A}$ onto  $R_{\restriction A'}$ which extends by the identity  on $F$ to a local isomorphism of $R$. We say that $R$ is \emph{almost-chainable}, resp. \emph{almost-monomorphic},  if there is some finite $F$ for which $R$ is $F$-chainable,  resp. $F$-monomorphic. The results mentionned above extend to these notions (see \cite {fraisse}).

\subsection{Monomorphic decomposition} \label{subsection:monomorphic}
Let $R$ be a relational structure.  A subset $V'$ of $V(R)$  is a \emph{monomorphic part} of $R$ if for every integer $k$ and every pair $A,~A'$ of $k$-element subsets of $V(R)$, the induced structures on $A$ and $A'$ are isomorphic whenever $A\setminus V'=A'\setminus V'$. A \emph{monomorphic decomposition} of $R$ is a partition $\mathscr P$ of $V(R)$ into monomorphic parts. A monomorphic part which is maximal for inclusion among monomorphic parts  is a \emph{monomorphic component} of $R$. The monomorphic components of $R$ form a  monomorphic decomposition of $R$ of which, every monomorphic decomposition of $R$ is a refinement (Proposition 2.12 of  \cite{P-T-2013}).  This decomposition is denoted by $\mathcal P(R)$ and called \emph{minimal} or \emph{canonical} and  
the number of its infinite blocks is called the \emph{monomorphic dimension} of $R$. We recall the following results, see Theorem 2.25 and Theorem 1.8 of \cite{P-T-2013}.
\begin{theorem}\label{thm:chainable} Every infinite component $A$ of a relational structure $R$ is $V(R) \setminus A$-chainable.
\end{theorem}

\begin{theorem}
 \label{thm:charcfinitemonomorph}
 A relational structure $R:=(V, (\rho_{i})_{i\in I})$ admits a finite
 monomorphic decomposition if and only if there exists a linear order $\leq$ on $V$ and a finite
 partition $(V_x)_{x\in X}$ of $V$ into intervals of $C:=(V, \leq)$ such that every local isomorphism of $C$ which preserves each interval is a local isomorphism of $R$. Moreover,
 there exists such a partition whose number of infinite blocks is the
 monomorphic dimension of $R$.
\end{theorem}

In other words, $R:=(V, (\rho_{i})_{i\in I})$ admits a finite monomorphic decomposition iff $R$ is freely  interpretable by a structure of the form $S:= (V, \leq, (u_1,  \dots, u_{\ell}))$ where $\leq$ is a linear order and $u_1, \dots, u_{\ell}$ are finitely many unary relations determining  a partition of $V$ into intervals of $(V, \leq)$.

A straightforward application of the Compactness Theorem of First Order Logic yields:

\begin{theorem} \label{compactness} A relational structure $R$ admits a finite monomorphic decomposition iff there is some integer $\ell$ such that every member of $\age(R)$ has a monomorphic decomposition into at most $\ell$ classes.
\end{theorem}

If $R$ is an ordered structure, it is obvious that an interval decomposition is a monomomorphic decomposition. The converse does not hold. However, the monomorphic components of $R$ which are infinite are intervals of the underlying chain (apply Theorem \ref{thm:chainable}). Hence, we have:

\begin{theorem} \label{thm:intervalsecomposition}An ordered relational structure $R$ has a finite interval decomposition iff it has a finite monomorphic decomposition; furthermore the monomorphic dimension of $R$ is the least number of infinite intervals of an interval decomposition.
\end{theorem}

The following property allows to test the existence of a finite monomorphic decomposition by looking at the   relations composing $R$, see Proposition 7.47 p.168  of \cite{oudrarthese} and \cite{oudrar-pouzet2014}.
\begin{proposition}
Let $m\in \mathbb N$. An ordered  relational structure $R:= (V, \leq, (\rho_i)_{i<m})$ has a finite monomorphic decomposition if and only if every $R_i:=(V, \leq, \rho_i)$ $(i<m)$ has such a decomposition.
\end{proposition}

\subsection{The role of well-quasi-order}\label{subsection:w.q.o.}

We recall that a poset $P$  is \textit{well-quasi-ordered (w.q.o.)} if $P$ contains no infinite antichain and no infinite descending chain. For more about this important notion, see \cite{fraisse}. Let
 $\mathscr{C}$ a subclass of  $\Omega_\mu$ and $P$ is a poset, we set $\mathscr{C}.P:=\{(R,f): R\in \mathscr{C},f: V(R)\rightarrow P\}$ and  $(R,f)\leq (R',f')$ if there is an embedding $h:R \rightarrow R'$ such that $f(x)\leq f'(h(x))$ for all $x\in V(R)$. We say that $\mathscr{C}$ is \textit{hereditary w.q.o.} if $\mathscr{C}.P$ is w.q.o.  for every w.q.o.  $P$.
 It is not known if $\mathscr{C}$ is  hereditary w.q.o. whenever  $\mathscr{C}. 2$  is w.q.o. for the $2$-element antichain $2$.  In fact, it is not even known if it suffices that  the collection of members of $\mathscr {C}$ with two distinguished constant is w.q.o. The preservation of wqo by labelling is central to the following result that we recall. If the signature $\mu$ is finite, a subclass of $\Omega_\mu$ which is hereditary and hereditary w.q.o. has finitely many bounds (\cite{pouzet 72}).

Let $\mathcal{P}:= (V_x)_{x\in X}$ be a partition of a set $V$ into disjoint sets. The \emph{shape} of   a finite subset $A$  of $V$ is the family $d(A):= (d_x(A))_{x\in X}$ where $d_x(A):= \vert A\cap V_x\vert$.

%
\begin{proposition}\label{prop-w.q.o.-ages}
The age of a relational structure admitting a finite monomorphic decomposition is hereditary w.q.o.
\end{proposition}
\begin{proof}
Let $R$ be such a relation, $V:= V(R)$, $\mathcal P:= (V_x)_{x\in X}$ be a monomorphic decomposition of $R$ into finitely many components, $\ell:= \vert X\vert$, $C_{x}$ be the chain $(\mathbb N, \leq)$ if $V_x$ is infinite and $C_{x}$ be the  finite chain $\{0, \dots, \vert V_x\vert\}$ otherwise; let  $C$ be the direct product  $\Pi_{x\in X} C_x$, ordered componentwise. Let $A$ and $A'$ be two finite subsets of $V$;  according to the definition of monomorphy, if they have  the same shape they are isomorphic. From which  follows that if $d(A)\leq d(A')$ then $R_{\restriction A} \leq R_{\restriction A'}$. Hence $\age (R)$ ordered by embeddability can be viewed as the surjective image of $C$. According to Dickson's lemma, $C$ is w.q.o., hence $\age(R)$ is w.q.o.   For the proof that $\age(R)$ is hereditary w.q.o.,  we use Theorem \ref {thm:charcfinitemonomorph} and Higman's theorem on words \cite{higman.1952}. This later result  asserts that on an  ordered alphabet $U$, the set $U^{<\omega}$ of finite words, ordered with the Higman ordering,  is w.q.o. provided that the alphabet is w.q.o. (we recall for two words $v:= (v_0, \dots, v_{n-1}), w:= (w_0, \dots, w_{m-1})$, the Higman ordering sets $v\leq w$ if there is an injective order-preserving map $h: [0, n[ \rightarrow [0, m[$ such that $v_i\leq w_{h(i)}$ for all $i<n$). Indeed, $R$ is freely interpretable by a structure of the form $S:= (V, \leq, (u_1,  \dots, u_{\ell}))$, where $\leq$ is a linear order and $u_1, \dots, u_{\ell}$  are finitely many unary relations forming a partition of $V$ into intervals of $(V, \leq)$. To prove that $\age(R). P$ is w.q.o. for every w.q.o. $P$, it suffices to prove that $\mathcal A(S). P$ is w.q.o. This conclusion follows from Higman's theorem. Indeed, let $L$ be the $\ell$-element antichain $\{1, \dots, {\ell}\}$ and $Q$ be the direct product $P\times L$. As  a product of two w.q.o.'s (in this case  a direct sum of $\ell$ copies of the w.q.o. $P$) $Q$ is w.q.o.. Hence, by Higman's theorem, the set $Q^{<\omega}$ of words on $Q$ is w.q.o.  To conclude, it suffices to observe that   $\mathcal A(S). P$ is the image  by an order-preserving map of a subset of $Q^{<\omega}$, and thus  it is w.q.o. This is easy: to a finite subset $A$ of $V$ and a map $f: A \rightarrow P$  we may associate the  word $w_A:= \alpha_0\cdots \alpha_{n-1}$ on $Q$ defined as follows: $n:= \vert A\vert$, $\alpha_i:= (f(a_{i}), \chi(a_i))$ for $i<n$, $\{a_0, \dots, a_{n-1}\}$ is the enumeration of $A$ w.r.t. the linear order $\leq$ on $V$ and $\chi(a):=j$ if  $a\in u_j$.
\end{proof}

By the main result of \cite{pouzet 72} mentionned at the begining of  this  subsection, we have:

\begin{corollary}\label{cor:w.q.o.-monomorphy}
If the signature is finite, the age of a structure with a finite monomorphic decomposition has  finitely many bounds.
\end{corollary}

We recall a  standard fact about ages, see the comment  made p.23, line 9 of \cite{P-T-2013} or  \cite{pouzet.81} 3.9, page 329.

\begin{theorem}\label {minimalage} If a hereditary classe $\mathscr C$ of finite relational structures with finite signature $\mu$ contains an   infinite antichain then it contains an age   which has infinitely many bounds in $\mathscr C$  and contains  no infinite antichain.
\end{theorem}
For reader's convenience, we repeat the proof.

\begin{proof}Consider the set $\mathfrak C$ of hereditary subclasses $\mathscr C'$  of $\mathscr C$, which have infinitely many bounds in $\mathscr C$; this means that, as a final segment of $\mathscr C$,  $\mathscr C\setminus \mathscr C'$ is not finitely generated. As it is easy to see, the intersection of any chain of members of $\mathfrak C$ belongs to $\mathfrak C$. According to Zorn's Lemma, $\mathfrak C$ contains a minimal element. If $\mathscr C_0$ is such an element, it has infinitely many bounds in $\mathscr C$. Furthermore, because it is minimal, $\mathscr C_0$ contains no infinite antichain. Finally, since the union of two hereditary classes with finitely many bounds has only finitely many bounds, $\mathscr C_0$ cannot be the union of two proper hereditary classes. This implies that  $\mathscr C_0$ is up-directed and thus is an ideal. According to Fra\"{\i}ss\'e's result already mentionned, this implies that $\mathscr C_0$ is the age of a relational structure.
\end{proof}

We extend the scope of Lemma \ref{lem:reduction} to monomorphic decompositions.
\begin{lemma} \label{lem:reductionbis}
If a hereditary class $\mathscr C$  of finite  relational structures with a finite signature $\mu$ contains for every integer $k$ some finite  structure which has no finite monomorphic decomposition  into at most $k+1$ blocks, then it contains a hereditary class  $\mathscr A$ with the same property which is minimal w.r.t. inclusion.
\end{lemma}

\begin{proof} If   $\mathscr C$ contains  no infinite antichain (with respect to embeddability) then the collection of its hereditary subclasses  is well-founded. And we are done. If $\mathscr C$ has an infinite antichain, Theorem \ref{minimalage} yields an   age $\mathcal A(R)$ with no infinite antichain and infinitely many bounds. Since $\mathcal A(R)$ has  infinitely many bounds, by Corollary \ref {cor:w.q.o.-monomorphy} $R$  cannot have a  finite monomorphic decomposition. Thus, by Theorem \ref{compactness}, $\mathcal A(R)$ has the same property as $\mathscr C$.  Since $\mathcal A(R)$ has no infinite antichain,  we are back to the first case.
\end{proof}

We may note that, as in the proof of Theorem \ref{minimalage}, each minimal class $\mathscr A$ in Lemma  \ref{lem:reductionbis} above is the age of a relational structure with no finite monomorphic decomposition.

A consequence is the following extension of Proposition \ref{prop-w.q.o.-ages}

\begin{proposition}\label{prop:hw.q.o.}
If the signature $\mu$ is finite and if  there is some integer $\ell$ such that every member of a hereditary class $\mathscr C$ has a monomorphic decomposition into at most $\ell$ classes then $\mathscr C$  is hereditary w.q.o.
\end{proposition}

\begin{proof}
According to Theorem \ref {minimalage} and Corollary \ref {cor:w.q.o.-monomorphy},  $\mathscr C$ contains no infinite antichain. This implies  that $\mathscr C$ is a finite union of up-directed initial segments. Each of these  initial   segments  is the age of a relational structure. By Theorem \ref{compactness} such a relational structure has a finite monomorphic decomposition. By Proposition \ref{prop-w.q.o.-ages} its age is hereditary w.q.o. Being a finite union of hereditary w.q.o. classes, $\mathscr C$ is hereditary w.q.o.
\end{proof}

\subsection{Profile}\label{subsection:profile}
The purpose of this subsection is to prove Theorem \ref{thm: polynomial-interval} and Theorem \ref{thm:minimal1}. To put these results in perspective, we recall few facts about the profile, see \cite{pouzet-profile2006} for a survey. First, the profile of an infinite age is non decreasing (Pouzet, 1971). A proof,
based on Ramsey's theorem,  appears in Fra\"{\i}sse's Cours de Logique  1971 \cite{fraisseclmt1} Exercice 8, p.113). In 1976 \cite{pouzet 76}  it was proved that  if $R$ is a relational structure, finite or not, on at least $2n+p$ elements, then $\varphi_{\age(\mathcal{R})}(n)\leq \varphi_{\age (\mathcal{R})}(n+p)$. The proof  uses  the non-degeneracy of the Kneser matrix of the $n$-element subsets of a $2n+p$ element set, a result  obtained independently by Gottlieb 1966 and Kantor 1972. In \cite{pouzet.81} was shown that the profile of a hereditary class union of infinite ages is non-decreasing. Relational structures whose profile is bounded,  hence  eventually constant,  were characterized as almost-chainable structures. Furthermore,  if the signature is finite and  bounded by $m$, then exists an integer  $s(m)$ such that if the profile of the age $\mathscr C$ of a relational structure of signature $\mu$ is constant on  $[n, . . . , n + s(m) + 1]$ then it is eventually constant. Hence, the growth of  $\varphi_{\mathscr C}$  is at least linear. More generally
the profile of the age $\mathscr C$ of a relational structure  with $\mu$ finite has either a polynomial growth or its  growth is faster than every polynomial \cite{pouzet.tr.1978} (the \emph{growth} of $\varphi_{\mathscr C}$ is {\it polynomial of degree $k$} if $a\leq \frac{\varphi_{\mathscr C}(n)}{n^k}\leq b$ where $a$ and $b$ are non-negative constants; it is  faster than every polynomial if  $\frac{\varphi_{\mathscr C}(n)}{n^k}$ goes to infinity with $n$ for every $k$).

If an infinite relational structure has a finite monomorphic decomposition into  $k+1$ blocks, then trivially, the  profile of $\mathcal A(R)$ is bounded by some polynomial whose degree is, at most, $k$. In fact, and this is the main result of \cite{P-T-2013},  this is eventually a quasi-polynomial whose degree is the monomorphic dimension of $R$ minus $1$. This is not eventually  a polynomial (for a simple example take   the direct sum of $k+1$ copies ($k\geq 1$) of the complete graph $K_{\infty}$ on infinitely many vertices). In the case of ordered structures, there is no symmetries: if two structures are isomorphic, the isomorphism is unique. We get easily the first part  of Theorem \ref{thm:poly-expo1} by using Lemma 2.15 of \cite{P-T-2013} valid for the notion of monomorphic decomposition that we record below in Lemma \ref{lemma.minimalGrowthRate}.

Let $(V_x)_{x\in X}$ be a partition of a set $V$. Given $d\in \mathbb N$,
call \emph{$d$-fat} a subset $A$ of $V$ such that, for all $x\in X$,
$d_x(A) \geq d$ whenever $A\not\supseteq V_x$.
\begin{lemma}
  \label{lemma.minimalGrowthRate}
  Let $R$ be an infinite relational structure on a set $V$ admitting a finite  monomorphic decomposition $(V_x)_{x\in X}$. Then, there exists some
  integer $d$ such that on every $d$-fat subset $A$ of $V$ (w.r.t.   $\mathcal P(R)$)  the partition $\mathcal P(R)_{\restriction A}$  induced by $\mathcal P(R)$ on $A$ coincides  with $\mathcal P(R_{\restriction A})$. In particular the shape of $A$ w.r.t. $\mathcal{P}(R)$ coincides with the shape of $A$ w.r.t. $\mathcal {P}(R_{\restriction A})$.
   \end{lemma}

There are very few cases for which  the profile $\varphi_{\mathscr C}$ of a hereditary  class $\mathscr C$ is polynomial (notably because  $\varphi_{\mathscr C}(0)=1$). We say that $\varphi_{\mathscr C}$ is \emph{eventually a polynomial} if for some non-negative integer $n$, the restriction of $\varphi_{\mathscr C}$ to the interval $[n, \rightarrow )$ of $\NN$ is a polynomial function.
This is the case if $\mathscr C$ is finite (the profile is eventually zero). More interestingly, this is the case if $\mathscr C$ is infinite but  each proper hereditary subclass of $\mathscr C$ is finite. In this case, we may observe that the profile is constant, equal to one. Indeed, via Ramsey's theorem, $\mathscr C$ is the age of a chainable relational structure. This simple observation allows to argue by induction in the proof of Theorem \ref{thm: polynomial-interval} below.

\begin{theorem} \label{thm: polynomial-interval}
If   there is some integer $k$ such that every member of a hereditary class $\mathscr C$ of ordered structures of finite signature $\mu$ has an interval decomposition into at most $k+1$ classes then the profile of $\mathscr C$  is eventuallly a polynomial of degree at most $k$.
\end{theorem}

\begin{proof}
According to Proposition \ref{prop:hw.q.o.},  $\mathscr C$ is hereditary w.q.o., hence w.q.o. Since the collection of hereditary subclasses is well founded we may  argue by induction. Thus,  we consider a hereditary class $\mathscr C'\subseteq \mathscr C$ such that $\varphi_{\mathscr {C''}}$ is eventually a polynomial for every proper hereditary subclass $\mathscr {C''}\subset  \mathscr C'$ and we prove that $\varphi_{\mathscr C'}$ is eventually  a polynomial. We may suppose that $\mathscr C'$ is non-empty (otherwise the  profile is undefined).

If $\mathscr C'$ is the union of two proper hereditary subclasses, say $\mathscr C _1'$ and $\mathscr C _2'$ then from the equality:
\begin{equation}\label{eq-profile}
\varphi_{\mathscr C'}=  \varphi_{\mathscr C_1'}+ \varphi_{\mathscr C'_2}-\varphi_{\mathscr C_1'\cap \mathscr C_2'}
\end{equation}
we get that $\varphi_{\mathscr C'}$ is eventually a polynomial. Thus we may suppose that $\mathscr C'$ is not the union of two proper hereditary subclasses.  Since $\mathscr C'$ is non-empty, this is the age of some relational structure $R$. According to Theorem \ref{compactness},  $R$ has a monomorphic decomposition   into finitely many blocks, hence by Theorem \ref{thm:intervalsecomposition} a finite interval decomposition. Let $\mathcal P:= (V_{x})_{x\in  X}$ be the canonical monomorphic decomposition, $X_{\infty}:= \{ x\in X:   V_{x}\;  \text{is infinite} \}$,  $K$ be the union of finite components of $R$,  $k+1:= \vert X_{\infty}\vert$ and $d$ given by Lemma \ref {lemma.minimalGrowthRate}. Let $\mathcal F$ be the set of $d$-fats finite subsets $A$ of $V$ containing $K$.  Due to Lemma \ref {lemma.minimalGrowthRate}, $\mathcal P(R_{\restriction A})$  is induced by $\mathcal P(R)$ and has the same number of classes as $\mathcal P(R)$. Since $R$ is ordered, two members $A, A' \in \mathcal F$ have the same shape iff  $R_{\restriction A}$ and $R_{\restriction A'}$ are isomorphic. Thus  the collection $\mathscr F:= \{R_{\restriction A}: AÂ \in \mathcal F\}$ is a non-empty final segment of $\mathscr C'$ and the map which associate  $(d_{x}(A)-d)_{x\in X_{\infty}}$ to $\mathcal R_{\restriction A}$ is an order isomorphism from $\mathscr F$ onto $\mathbb N^{k+1}$.  Let $\vartheta$ be the generating function  of $\mathbb N^{k+1}$ (which counts for each integer $m$ how many elements of $\mathbb N^{k+1}$ have level  $m$; we have $\vartheta (m)={ m+k \choose k}$ and $\varphi_{\mathscr F}(n)= \vartheta (n-(k+1)d- \vert K\vert)$;  hence $\varphi_{\mathscr F}$ is a polynomial of degree $k$. Since  $\mathscr C'\setminus \mathscr F$ is  a proper hereditary class of $\mathscr C'$,  induction asserts that   its profile is eventually a polynomial. Hence the profile of $\mathscr C'$ is eventually a polynomial.
\end{proof}

\medskip


Corollary \ref{cor:minimal} can be extended as follows.

\begin{theorem}\label{thm:minimal1}
If the profile of a hereditary class $\mathscr C$ of finite relational structures (with a finite signature $\mu$) is not bounded above by a polynomial then it contains a hereditary class  $\mathscr A$ with this property which is minimal w.r.t. inclusion. Hence, this class is an age.
\end{theorem}
The proof follows the same line as the proof of Lemma \ref{lem:reductionbis}. The key ingredient replacing
Corollary \ref{cor:w.q.o.-monomorphy} is the following lemma:

\begin{lemma}
If the signature is finite, the age of a structure with polynomially bounded profile has  finitely many bounds.
\end{lemma}

The proof of this fact is far from being trivial. See \cite{P-T-2013} for the  main  ingredients.


\section{Monomorphic decomposition, multichainability and invariant structures}

In this section, we present an alternative definition of monomorphic decomposition, simpler to use than the original one. We recall the definition of almost-multichainable structures and give a test for those with  no finite monomorphic decomposition. We conclude with a short presentation of  the notion of invariant structure.

 \subsection{Monomorphic decomposition via an equivalence  relation}\label{subsection:canonical}
 We   revisit the notion of monomorphic decomposition. We  define it in a direct way, via an equivalence relation. For more, see Section 7.2.5 and Chapter 8 of \cite{oudrarthese}.

Let $R$ be a relational structure and $V:= V(R)$. Let $x$ and  $y$ be two elements of $V$. Let  $A$ be a finite subset of $V \setminus \{x, y\}$, we say that $x$ and  $y$  are \emph{$A$-equivalent} and we  set  $x \simeq_{A, R} y$ if   the restrictions of $R$ to   $\{x\}\cup A$  and $\{y\}\cup A$ are isomorphic. Let  $k$ be a non-negative integer, we set $x \simeq_{k, R} y$ if $x \simeq_{A, R} y$ for every $k$-element  subset $A$ of   $V(R)\setminus \{x,y\}$.  We set $x\simeq_{\leq k, R}y$ if $x\simeq_{k', R}y$ for every $k'\leq k$ and $x\simeq_{ R}y$ if $x \simeq_{k, R} y$ for every $k$.

For example, the 0-equivalence relation is an equivalence relation, two vertices $x$ and $y$ being
equivalent if and only if the restrictions of $R$ to $x$ and $y$ are isomorphic.

\begin{lemma}\label{lem:equivalence}
The relations $\simeq_{k, R}$,  $\simeq_{\leq k, R}$ and $\simeq_{ R}$ are equivalence relations on $V$. Furthermore, each equivalence class of $\simeq _R$ is a monomorphic part of $R$ and each monomorphic part is included into some equivalence class.
\end{lemma}

\begin{proof}
We prove that $\simeq_{k, R}$ is an equivalence relation on $V$. For that it suffices
to check that it is transitive. Let $ x,y, z\in V$ with $x\simeq_{k, R} y$ and $y\simeq_{k, R} z$.
We check that $x\simeq_{k, R} z$. We may suppose these elements pairwise distinct.
Let $A$ be a $k$-element subset of $V\setminus \{x,z\}$. Case 1) $y\not \in A$. In this case
$R_{\restriction A\cup\{x\}}$ is isomorphic to $R_{\restriction A\cup\{y\}}$ and $R_{\restriction A\cup\{y\}}$ is isomorphic to $R_{\restriction A\cup\{z\}}$.
Hence $R_{\restriction A\cup\{x\}}$ is isomorphic to $R_{\restriction A\cup\{z\}}$; and thus $x\simeq_{k, R} z$. Case 2) $y\in A$.
Set $A':= (A\setminus \{y\})\cup \{z\}$. Since $x\simeq_{k, R} y$, we have $R_{\restriction A'\cup\{x\}}$ isomorphic to
$R_{\restriction A'\cup\{y\}}= R_{\restriction A\cup \{z\}}$. Similarly, setting $A'':= (A\setminus \{y\})\cup \{x\}$ then, since
$y\simeq_{k, R} z$,  we have $ R_{\restriction A\cup \{x\}}=R_{\restriction A''\cup\{y\}}$ isomorphic to $R_{\restriction A''\cup\{z\}}$. Since
$A'\cup\{x\}=A''\cup\{z\}$ the conclusion follows. The fact  that  $\simeq_{\leq k, R}$ and $\simeq_{ R}$ are also equivalence relations follows.

Let $C$ be an equivalence class of $\simeq_{R}$. We prove that $C$ is a monomorphic part.  Let   $A, A'\in [V]^h$  such that $A\setminus C=A'\setminus C$. Our aim is to prove that $R_{\restriction A}$ and $R_{\restriction A'}$  are isomorphic. Let $\ell := \vert A\setminus A'\vert$. If  $\ell=0$, $A=A'$, there is nothing to prove. If $\ell =1$,  then $A= \{x\}\cup (A\cap A')$ and $A'=\{y\}\cup (A\cap A')$, with $x,y\in C$; in this case $R_{\restriction A}$ and $R_{\restriction A'}$  are isomorphic since $x\simeq_{R} y$. If $\ell >1$, set $K:=A \setminus C$ and $k:= \vert K\vert$. We may find a sequence of $(h-k)$-element subsets of $C$, say $A_0, \dots A_i, \dots, A_{\ell-1}$ such that $A_0= A\cap C$,  $A_{\ell-1}= A'\cap C$ and the symmetric difference of $A_i$ and $A_{i+1}$ is $2$. 
We have $ R_{\restriction  A_i\cup K}=  R_{\restriction A_{i+1}\cup K}$ for $i<\ell-1$. Hence $R_{\restriction A}$ and $R_{\restriction A'}$ are isomorphic as required. Since the equivalence  classes are  monomorphic parts  they form a monomorphic decomposition. Trivially, the elements of a monomorphic part  are pairwise equivalent w.r.t. $\simeq_{R}$ hence contained into an equivalence class.
\end{proof}

From Lemma \ref{lem:equivalence} we obtain:

\begin{proposition}\label{prop: equiv=decomp}
The equivalence classes of $\simeq_{R}$ form a decomposition of $R$ into monomorphic parts and every decomposition into monomorphic parts is finer. Thus the decomposition of $R$ into equivalence classes of $\simeq _{R}$
coincides with the decomposition of $R$ into monomorphic
components.
\end{proposition}

Some properties mentionned in Subsection \ref{subsection:monomorphic}, e.g., Theorem \ref{compactness} and Lemma \ref{lemma.minimalGrowthRate},   follow easily from  Proposition \ref{prop: equiv=decomp}.

 \subsection{Hypomorphy  and equivalence}\label{subsection: hypomorphy}
Let $R$ and $R'$ be two relational structures on the same set $V$ and $k$ be an integer. We say, with Fra\"{\i}ss\'e and Lopez \cite{fraisse-lopez},  that $R$ and $R'$ are \emph{$k$-hypomorphic} if the restrictions $R_{\restriction A}$ and $R'_{\restriction A}$ are isomorphic for every $k$-element subset $A$ of $V$. Let $R$ be a relational structure  on $V$ and $x,y\in V$. Set $V':= V\setminus \{x,y\}$. By identifying $x$ and $y$ to some element $z$ (eg $z:= \{x,y\}$) we get two structures $R(x)$ and $R(y)$ on $V' \cup \{z\}$. Formally,  $R(x)$ is such that the map which fixes $V'$ pointwise and send $x$ on $z$ is an isomorphism of $R_{\restriction  {V\setminus \{y\}}}$ onto $R(x)$; similarly for $R_{\restriction {V\setminus \{x\}}}$ and $R(y)$.

Since $R(x)$ and $R(y)$ coincide on $V'$ we have immediately:

\begin{lemma}\label{lem:iteration}Let $x,y\in  V$ and $k\in \mathbb N$ then  $x\simeq _{k, R} y$ iff $R(x)$ and $R(y)$ are $(k+1)$-hypomorphic. \end{lemma}

It was shown that \emph{two relational structures which are $(k+1)$-hypomorphic  are $k$-hypomorphic provided that their domain has at least $2k+1$ elements}; see \cite{pouzet 79} Corollary 2.3.2 and also \cite{pouzet 76} where this fact is obtained  from a property of  incidence matrices  due to Gottlieb \cite{gottlieb} and Kantor \cite{kantor}.
 With this property and Lemma \ref{lem:iteration},  we obtain:

\begin{corollary}\label{lem:gottlieb-kantor}
The equivalence relations $\simeq_{k, R}$ and $\simeq_{\leq k, R}$ coincide provided that $\vert V(R)\vert \geq 2k+1$.
 \end{corollary}

The notion of hypomorphy has been fairly well studied, particularly for binary structures. Let us say that two structures $R$ and $R'$ on the same set $V$ are \emph{hypomorphic} if they are $k$-hypomorphic for every $k< \vert V\vert$ and that $R$ is \emph{reconstructible} if $R$ is isomorphic to every $R'$ which is hypomorphic to $R$.
Lopez  showed that every finite binary structure on at least $7$ elements  is  reconstructible (\cite{lopez 72} 1972 for binary relations, \cite{lopez 78} 1978 for binary structures, see the exposition in  \cite{fraisse-lopez}). From Lopez's result     and Lemma \ref{lem:iteration}  follows:

\begin{theorem}\label{thm:binarytreshold}
On  a binary structure $R$, the equivalence relations $\simeq _{\leq 6, R}$ and $\simeq _{R}$ coincide.
\end{theorem}

  Ille (1992) \cite {ille 92} showed  that for every integer $m$ exists some integer $s(m)$ such that all finite ordered structures whose  maximum of the arity is $m$ and  cardinality at least  $s(m)$  are reconstructible. Thus,  similarly to Theorem \ref{thm:binarytreshold}, we  have:

  \begin{theorem}\label{thm:marytreshold} For every integer $m$ there is an integer $i(m)$, $i(m)\leq s(m)$,  such that the equivalence relation $\simeq _{\leq i(m),R}$ and $\simeq _{R}$ coincide on every  ordered structure $R$  of arity at most $m$ and cardinality at least $i(m)$.
\end{theorem}

The threshold for binary structures is examined in Chapter 8 of \cite{oudrarthese}. It is shown that the value $6$ in Theorem \ref {thm:binarytreshold}  can be  replaced by $3$ for binary structures and by $2$ for ordered binary structures (a result obtained also by Y.Bouddabous \cite{boudabbous}).

 The second author  showed in 1979 \cite {pouzet 79} that there is no reconstruction threshold  for ternary relations. He constructed for each integer $m$ a pair of hypomorphic ternary relations $S$ and $S'$ on a set $W$ of size $2m$  which are not isomorphic. In his example,  $W$ contains some  element $z$ (in fact two) such that $S$ and $S'$ coincide on $W\setminus \{z\}$. Hence we may duplicate $z$ and built a ternary relation $R$ on $2m+1$ elements having two elements $x,y$ such that $x\not \simeq_Ry$ but $x\simeq_{\leq 2m-1, R}y$.

 We do not know if Theorem \ref{thm:binarytreshold}  extends to infinite ternary structures. We can only offer two opposite problems.

\begin{problem} Existence of a relational structure  $R$ with  a finite signature such that $\simeq_R$ has infinitely many equivalence classes and $\simeq_{k, R}$ has only finitely many  equivalence classes for each integer $k$.
 \end{problem}
 On the opposite:

 \begin{problem}\label{problem:threshold} Existence for each integer $m$ of an integer $o(m)$ such that every infinite class of $\simeq_{o(m),R}$-equivalence of a relational structure with arity at most  $m$ is a class of $\simeq_R$-equivalence.
 \end{problem}

Then:
 \begin{lemma}Under a positive solution to Problem \ref{problem:threshold}, every relational structure $R$ of arity at most $m$ has the same infinite classes of $\simeq_R$-equivalence and of $\simeq_{k, R}$-equivalence for each $k\geq o(m)$. As a consequence, if $R$ has  infinitely many classes of $\simeq_R$-equivalence, it has infinitely many classes of $\simeq_{k,R}$-equivalence.
 \end{lemma}

\begin{proof} Since $V(R)$ is infinite, a class of $\simeq_{k, R}$-equivalence  is a class of $\simeq_{k', R}$-equivalence for every $k'\leq k$ by Corollary \ref{lem:gottlieb-kantor}. In particular, a class $A$ of $\simeq_{k, R}$-equivalence  is a $\simeq_{o(m), R }$-equivalence class.  Assuming a  positive result to the problem above, if $A$  is  infinite this  is    a class of $\simeq_R$-equivalence. Conversely,  if $A$ is a class of  $\simeq_R$-equivalence,  it is contained in a class $A'$ of $\simeq_{o(m), R}$-equivalence. If $A$ is infinite, then $A'$ is infinite, hence this is a class of  $\simeq_R$-equivalence by the positive result of Problem \ref{problem:threshold}, thus $A'=A$. The consequence is immediate.
\end{proof}

\begin{proposition}\label{prop:tresholdwqo}
If $\mathcal A(R)$, the age of a relational structure $R$,  is hereditary w.q.o.  then  the equivalences relations $\simeq_R$ and $\simeq_{\leq k, R}$ coincide for some nonnegative integer $k$. \end{proposition}

\begin{proof}Let $V:= V(R)$. Suppose that for every nonnegative integer $k$ there are $x_k, y_k \in V$ such that $x_k \simeq_{k, R} y_k$ and  $x_k \not \simeq_{R} y_k$; let $A_k\subseteq V\setminus \{x_k,y_k\}$ witnessing  the fact that $x_k \not \simeq_{R} y_k$. For $k$, set  $S_k:= (R_{\restriction A_k \cup \{x_k, y_k\}}, x_k, y_k)$ be the relational structure made of the restriction of $R$ to  $A_k \cup \{x_k, y_k\}$
and the  two distinguished element $x_k$, $y_k$. Suppose that   for some $k$ and $k'$, $S_k$ is embeddable into $S_{k'}$. We claim that $k'\leq \vert A_k\vert$.  Indeed, suppose the contrary. Let $f$ be such an embedding and $A'_k:= f(A_k)$. Since  $x_{k'}\simeq_{k', R} y_{k'}$, and  $k'> \vert A'_k\vert $ then, $R_{\restriction A'_{k}\cup \{x_{k'}\}}$ and $R_{\restriction A'_k\cup \{y_{k'}\}}$ are isomorphic too. But then $R_{\restriction A_k\cup \{x_k\}}$ and $R_{\restriction A_k\cup \{y_k\}}$ are isomorphic. This  contradicts the fact that
$A_k$ witnesses  that $x_k \not \simeq_{R} y_k$. This proves our claim. Now, if $\mathcal A(R)$ is hereditarily w.q.o., in particular, if the collection of members of $\mathcal R$ with two distinghuished elements is w.q.o, then from the sequence $(S_k)_{k\in \NN}$ we may extract an increasing subsequence. This fact violates our claim. The conclusion of the proposition follows.
\end{proof}

\subsection{Almost-multichainability}\label{subsection:almostmulti}

 An  \emph{almost-multichain} is a relational structure $S:= (V, \leq, \rho_L, U_1,\dots, U_k,  a_1, \dots a_{\ell})$, where $V$ is the disjoint union of a finite $\ell$-element set $F$ and a set  of the form $L\times K$ where  $K$ is a $k$-element set $\{v_1, \dots, v_k\}$;  the relation $\rho_L$ is an equivalence  relations on $V$ such that all the elements of $F$ are inequivalent and inequivalent to those of $L\times K$, furthermore two elements   $(x,y)$ and $(x',y')$ of $L\times K$ are equivalent w.r.t. $\rho_L$ if $x=x'$; the relations $U_1, \dots, U_k$ are unary relations, with  $U_i(x, y)=1$ if $y=v_i$,  the relations $a_1, \dots, a_\ell$ are constants  defining the $\ell$-element set  $F$; the relation $\leq $ is an order on $V$ for which all the elements of $F$ are minimal and below all elements of $L\times K$, all blocks of $\rho_L$ are antichains for the ordering and the family of these blocks is totally ordered by $\leq$. A relational structure $R$ freely interpretable by  an almost-multichain $S$ is \emph{almost-multichainable}.

 The  notion of almost-multichain, a bit complicated, allows to use Compactness Theorem of First Order Logic.  The relational structure $R$ is freely interpretable by $S$  just  in case for every  local isomorphism $h$ of the chain $C:= (L, \leq)$ (where $\leq $ denotes the order on $L$ induced by the order on $V$),  the map $(h, 1_K)$
extended by the identity on $F$ is a local isomorphism of $R$ (the map
$(h, 1_K)$ is defined by $(h, 1_K)(x, y):= (h(x), y)$).

Let $m$ be an  integer. We will say that the multichain $S$ and every $R$  freely  interpretable by $S$
\emph{has type at most}  $(\omega, m)$  if $C$ has order type $\omega$ and $\vert F\vert+ \vert K\vert\leq m$.

%
The notion of almost-multichainability was introduced in \cite{pouzet.tr.1978} (see \cite{pouzet-profile2006} for further references and discussions).  The special case $\vert K\vert =1$ is
the notion of almost-chainability introduced by Fra\"{\i}ss\'e. It seems to be a
little bit hard to swallow, still it is the key in proving the second part of Theorem \ref{thm:poly-expo1} as well that Theorem \ref{theo:base}.

We recall the following result (Theorem 4.19, p.265,  of \cite{pouzet-profile2006}); we just give a hint.

\begin{proposition}\label{prop:hereditwqo} The age $\mathcal A(R)$ of an almost-multichainable relation is hereditary w.q.o.
\end{proposition}

\begin{proof}
Let $S:= (V, \leq, \rho_L, U_1,\dots, U_k,  a_1, \dots a_{\ell})$ be an almost-multichain interpreting $R$. Let $P$ be w.q.o. In order to prove that $\mathcal A(R). P$ is w.q.o. it suffices to prove that $\mathcal A(S). P$ is w.q.o. For
that, let $A$ be the alphabet  whose letters are the non-empty subsets of $K$ (that is
$A := \powerset (K) \setminus \{\emptyset\})$. We order $A$ by inclusion and
extends this order to  $A^{*}$, the set of words over $A$, via the Higman's ordering.  Let $F:= \{a_1, \dots a_{\ell}\}$. Order $\powerset (F)$  by inclusion. Then $\mathcal A(S)$ is isomorphic to the direct product $\powerset (F) \times A^{*}$. Via Higman's theorem  on words, $\mathcal A(S)$ is w.q.o.  Then, replace  $A$ by the set of $f: X\rightarrow  P$ such that $X\in A$ and  associate corresponding  set of words. An other application of Higman's theorem  yields that $\mathcal A(S). P$ is  w.q.o.
\end{proof}

 Applying  Proposition \ref {prop:tresholdwqo}, Proposition \ref{prop:hereditwqo} and  the main result of \cite{pouzet 72}, we get:

\begin{proposition}\label{prop:treshold-bounds} If $R$ is an almost-multichainable relational structure then  there is some integer $k$ such that $\simeq_R$ and $\simeq_{\leq k, R}$ coincide. If moreover  the arity is finite, then $R$ has only finitely many bounds.
\end{proposition}

An almost-multichainable relation may have a finite monomorphic decomposition (e.g., the direct sum of finitely many copies of the countable complete graph), or no finite monomorphic decomposition (e.g.,  the direct sum of infinitely many copies of a complete graph on a finite set of vertices). Here is a  simple fact.

\begin{lemma}\label{lem:test1}  Suppose that $R$ is a relational structure freely interpretable by an almost  multichain  $S:= (V, \leq, \rho_L, U_1,\dots, U_k,  a_1, \dots a_{\ell})$. Then,
for every $y\in K$, either all elements of $L\times \{y\}$ are  $\simeq_R$-equivalent, or there is some interval $L'$ of the chain $L$ such that $L\setminus L'$ is finite and all the elements of $L'\times \{y\}$ are $\simeq_R$-inequivalent. 
\end{lemma}

\begin{proof}
Let $F:= \{a_1, \dots, a_{\ell}\}$. Let $k$ be given by Proposition \ref{prop:treshold-bounds}.

 \begin{claim}\label{claim:treshold-bounds}
Suppose that there  are  $u< v\in L$ such that in the chain $L$,  the intervals $(\leftarrow u[$ and $]v \rightarrow)$ have at least $k$ elements and that $(u,y ) \simeq_R  (v,y)$.  Then,  all the elements of $L\times \{y\}$ are  $\simeq_R$-equivalent.
\end{claim}
\noindent{\bf Proof of Claim  \ref{claim:treshold-bounds}.} Being infinite, $L$ contains either a copy of the chain $\NN$ or a copy of its dual. We may suppose that it contains a copy of $\NN$ (otherwise replace $L$ by its dual). Let $L':= x_0<\dots <x_n<x_{n+1}< \dots$ be such a chain. Let $R':=R_{\restriction V'}$ where $V':= F\cup (L'\times K)$. Let $n\in \NN$. Then $(x_n, y) \simeq_{k, R'} (x_{n+1},y)$. Otherwise, there is some subset $A' \subset V'\setminus \{(x_n, y), (x_{n+1},y)\}$ with at most $k$ elements such that $R'_{\restriction \{(x_n,y) \}\cup A'}$ and $R'_{\restriction \{(x_{n+1} ,y)\}\cup A'}$ are not isomorphic. But, since $A'$ is either on the left of $\{x_{n}\} \times K$ or on the right of  $\{x_{n+1}\} \times K$, there is a local isomorphism of $R$ defined on $A'\cup \{(x_n, y), (x_{n+1},y)\}$ mapping $(x_n, y)$ on $(u, y)$ and $(x_{n+1},y)$ on $(v,y)$. But then  $(u,y ) \not \simeq_R  (v,y)$, contradicting our supposition.  Now, from the transitivity of  $\simeq_{\leq k,  R'}$ we get that
all the elements of $L'$ are $\simeq_{\leq k,  R'}$-equivalent.  We deduce from that fact that   all elements of $L\times \{y\}$ are  $\simeq_R$-equivalent. If not,  there are  $u, v\in L$ such that  $(u,y) \not \simeq_R (v,y)$. Let   $A$ be a finite subset of $V\setminus \{(u,y), (v,y)\}$  witnessing that fact. Since $\simeq_{R}$ and $\simeq_{\leq k,R}$ coincide, we may suppose that $A$ has at most $k$ elements. Since $R$ is almost-multichainable and $L'$ is an infinite subchain of  $R$, $R$ and $R'$ have the same age; in fact there is a local isomorphism of $R$ in $R'$ carrying  $A\cup  \{(u,y), (v,y)\} $ on $A'\cup  \{(u',y), (v',y)\}$ with $u',v'\in L'$. Since, we may choose $u'=x_n$ and $v'=x_{n+1}$ for some $n$, we get  $(u',y)\not \simeq_R (v',y)$. A contradiction. This proves our claim.
The proof of the lemma follows immediately from the claim.  \end{proof}
%

\begin{question} Is it true that all elements of  $L\times \{y\}$ are either $\simeq_{R}$-inequivalent or all are $\simeq_{R}$-equivalent?
\end{question}

As defined in the introduction, let $\mathscr S_{\mu}$ be the class of all relational structures of signature $\mu$, $\mu$ finite, without any finite monomorphic decomposition.

\begin{corollary}\label{lem:test}  Let $R$ be an infinite almost-multichainable structure   on  $V:= F\cup (L\times K)$. Then $R\in \mathscr S_{\mu}$ iff there are $x,x'\in L$ and $y\in K$ such that $(x,y)\not \simeq_{R} (x',y)$. This is particularly the case if for $A_x:= F\cup (\{x\}\times (K \setminus \{y\})$, the restriction of $R$ to $\{(x,y)\}\cup A_{x}$ and to $\{(x',y)\}\cup A_{x}$, are not isomorphic.
\end{corollary}


\begin{proposition}\label{prop:number}
If the signature $\mu$ is finite, there are finitely many pairwise non isomorphic almost-multichainable  structures of signature  $\mu$  with type at most $(\omega, m)$.
\end{proposition}
\begin{proof} Let $\ell$ be the maximum of $\mu$. Let $m$ be a non-negative integer, let $V:= F\cup (L\times K)$ where $F$ and  $K$ are two finite sets such that $\vert F\vert +\vert K\vert \leq m$ and $\leq$ be a linear order  on $L$ of type $\omega$. Let $L'$ be a  $\ell$-element subset of $L$.
Observe that if $R$ and $R'$  are two almost-multichainable  relational structures on $V$ with this decomposition which coincide on $F\cup (L'\times K)$ they are equal. The number of relational structures of signature $\mu$ on a finite set is finite, hence the number of those $R$ must be finite.\end{proof}

\subsection{Invariant structures}

A tool  to construct almost-multichainable structures is Ramsey's theorem. We will use that tool  via
the notion of invariant structure.  We refer to \cite{charretton-pouzet} for the notion and to  \cite{Bou-Pouz} and \cite{oudrar} for   examples of use.

Let  $C:=(L,\leq )$ be  a chain.
For each non-negative integer $n$, let $[C]^{n}$ be the set of $n$-tuples $\overrightarrow{a}:=(a_{1},...,a_{n})\in L^{n}$ such that $a_{1}<...<a_{n}$. This set will be
identified with the set of the $n$-element subsets of $L$. For every local automorphism $h$ of $C$ with domain $D$, set $h(\overrightarrow{a}):=(h(a_{1}),...,h(a_{n}))$ for every $\overrightarrow{a}\in \lbrack D]^{n}$.

Let $\mathfrak{L}:=\left\langle C,{R},\Phi \right\rangle $ be
a triple made of a chain $C$ on $L$,  a relational structure ${R}:=(V,(\rho _{i})_{i\in I})$ and a set $\Phi $ of maps, each one being a map
$\psi $ from $[C]^{a(\psi )}$ into $V$, where $a(\psi )$ is an integer, the \emph{arity} of $\psi$.

We say that
$\mathfrak{L}$ is \emph{invariant} if:
\begin{equation}\label{eq:invariance}
\rho _{i}(\psi _{1}(\overrightarrow{\alpha }_{1}),...,\psi _{m_{i}}(%
\overrightarrow{\alpha }_{m_{i}}))=\rho _{i}(\psi _{1}(h(\overrightarrow{%
\alpha }_{1})),...,\psi _{m_{i}}(h(\overrightarrow{\alpha }_{m_{i}})))
\end{equation}

\begin{equation}\label{eq:invarianceequality}
\psi_1(\overrightarrow{\alpha_{1}})=\psi_2 (\overrightarrow{\alpha_{2}})\;  \text {iff} \;
\psi_1((h(\overrightarrow{\alpha}_{1}))=\psi _{2}(h(\overrightarrow{\alpha}_{2}))
\end{equation}

for every $i\in I$ and every local automorphism $h$ of $C$ whose domain contains $\overrightarrow{\alpha }_{1},...,\overrightarrow{\alpha }_{m_{i}},$
where $m_{i}$ is the arity of $\rho _{i}$, $\psi _{1},...,\psi _{m_{i}}\in\Phi ,$ $\overrightarrow{\alpha }_{j}\in \lbrack C]^{a(\psi _{j})}$ for $j=1,...,m_{i}.$

These conditions expresse  the fact that each $\rho _{i}$ is invariant under the transformation of the $m_{i}$-tuples of $V$ which are induced on $V$ by the local automorphisms of $C$.  The second condition is satisfied if some relation $\rho_i$ is the equality relation or an order. For example, if ${\rho}$ is a binary relation and $\Phi =\{\psi\} $ then
\begin{equation*}
\rho(\psi (\overrightarrow{\alpha }),\psi (\overrightarrow{\beta }))=%
\rho(\psi (h(\overrightarrow{\alpha })),\psi (h(\overrightarrow{\beta})))
\end{equation*}%
means that, $\rho (\psi (\overrightarrow{\alpha }),\psi (\overrightarrow{\beta }))$ depends only upon the relative positions of $\overrightarrow{\alpha }$ and $\overrightarrow{\beta }$ on the chain $C.$

If $\mathfrak{L}:=\left\langle C,{R},\Phi \right\rangle $ and $A$ is a subset of $L$, set $\Phi _{\upharpoonleft _{A}}:=\{\psi _{\upharpoonleft _{\lbrack A]^{a(\psi )}}}:\psi \in\Phi \}$ and $\mathfrak{L}_{\upharpoonleft _{A}}:=\left\langle
C_{\upharpoonleft _{A}},{R},\Phi _{\upharpoonleft_{A}}\right\rangle $ the restriction of $\mathfrak{L}$ to $A.$

As we are going to see, almost-multichainable structures fit in the frame of invariant structures.

Let  $S:= (V, \leq, \rho_L, U_1,\dots, U_k,  a_1, \dots a_{\ell})$ be an almost-multichain as defined in subsection \ref{subsection:almostmulti}. Let $C:= (L, \leq)$. For each $k\in
K$,  let $f_k$ be the unary map from $C$ to $V$ defined by $f_k(a):= (a, k)$; also  for each $a_j$, $1\leq j\leq\ell$,  let $f_j$ be the $0$-ary map from $C^{0}$ taking the value $a_j$. Let $\Phi:= \{f_i: i\in K\cup \{1, \dots, \ell\}\}$.

As it is easy to check, we have:

\begin{proposition} A relational structure  $R$ on $V$ is freely interpretable  by $S$, hence almost-multichainable,   if and only if $\mathfrak{L}:=\left\langle C,{R},\Phi \right\rangle $  is invariant.
\end{proposition}

The following result \cite{charretton-pouzet}, consequence of Ramsey's theorem, will  be applied in the next section.

\begin{theorem}\label{thm:ramsey-invariant}
Let $\mathfrak{L}:=\left\langle C,{R},\Phi \right\rangle $ be a structure such that the domain $L$ of $C$ is infinite, ${R}$
consists of finitely many relations and $\Phi $ is finite. Then there is an infinite subset $L^{\prime }$ of $L$ such that $\mathfrak{L}_{\upharpoonleft
_{L^{\prime }}}$ is invariant.
\end{theorem}

\section{Proofs of Theorems \ref {theo:base} and \ref{theo:expo}}\label{section:proofs}

 Let $k$ be an integer. Let $\mathscr S_{\mu, k}$ be the class of structures $R$ of signature $\mu$ such that $\simeq_{\leq k, R}$ has infinitely many classes.

\begin{theorem}\label{mainthm}
Let $\mu$ be a finite signature and $k$ be a nonnegative  integer, then $\mathscr S_{\mu, k}$ contains a finite subset $\mathfrak A$ made of almost multichainable structures  such that every member of $\mathscr S_{\mu, k}$ embeds some member of $\mathfrak A$.
\end{theorem}

%

\begin{proof} We prove that every member of $\mathscr S_{\mu, k}$ embeds some member of $\mathscr S_{\mu,k}$ which is almost-multichainable and has type at most $(\omega, k+1)$. This suffices since according to  Proposition \ref{prop:number} the number of these structures is finite.

Let  $R\in \mathscr S_{\mu, k}$. Let $V$ be the domain of $R$. Then $\simeq_{\leq k, R}$ has infinitely many classes.  Pick  an infinite sequence $(x_{p})_{p\in \mathbb N}$ of pairwise inequivalent elements of $V$. For each pair $(p,q)$, $p<q$,  we can
find a subset $\mathcal F(p,q)$ of $V\setminus \{x_p, x_q\}$ with  at most  $k$ elements witnessing the fact that $x_{p}$ and $x_{q}$ are inequivalent. Hence,  we may find a family $ \Phi$ of maps $f:\mathbb{N}\rightarrow V$,  $g_{i}:[\mathbb{N]}^{2}\rightarrow V$ for $i=1,\dots, k$ such that for each $p<q\in \mathbb{N}$, $x_{p}=f(p)$, $\mathcal F(p,q)=\{g_{i}(p,q),~i=1,\dots, k\}$.

Let $C:= (\mathbb N, \leq)$ and $\mathfrak {L}:=\left\langle C,{R},\Phi \right\rangle $.  Ramsey's theorem under the form of Theorem \ref{thm:ramsey-invariant} asserts that there is some infinite subset $X$
of $\mathbb{N}$ such that $\mathfrak {L}_{\restriction X}$ is  invariant.

Without loss of generality, we may suppose that
$X= \mathbb N$, that is, $\mathfrak {L}$ is invariant. Furthermore, we may suppose that all members of $C$ are distinct. With these assumptions, we built  some  restriction of $R$ which is almost-multichainable and has type at most $(\omega, k+1)$.

The proof is based on three claims below.

\noindent\begin{claim}\label{basicclaim}
Suppose that  $\mathfrak {L}$ is  invariant.
\begin{enumerate}
\item Let $i\in \{1, \dots, k \}$, then $g_i(p, q)\not = f(r)$  for  some $p<q$ and $r$ distinct from $p, q$.
\item Let $i,j\in \{1, \dots, k\}$. Then $g_i=g_j$ iff there are $p<q$ such that $g_i(p,q)=g_j(p,q)$.
\item Let $i\in \{1,\dots, k\}$. Then $g_i$ is constant iff there are  integers $p<q$ and $p'<q'$ with $p\not =p'$ and $q\not = q'$ such that $g_i(p,q)=g_i(p',q')$.
\item If the restrictions $R_{\restriction \{f(p) \cup \mathcal F(p,q)\}}$ and  $R_{\restriction \{f(r) \cup \mathcal F(p,q)\}}$ are isomorphic for some integers $p<q<r$ then the corresponding restrictions for $p'<q'<r'$ are isomorphic.
\end{enumerate}
\end{claim}
\noindent{\bf Proof of Claim \ref{basicclaim}}. $(1)$ Suppose  $g_i(p, q) = f(r)$  for some  $p<q$ and $r$ distinct from $p$ and $q$. Then, we may find $p',q', r', r''$ and two  local isomorphisms defined on  $p,q,r$ sending these elements on $p', q',r'$ and on $p',q', r''$ respectively. Fom Equation (\ref{eq:invarianceequality}) defining the invariance, we get $g_i(p', q') = f(r')$ and $g_i(p', q') = f(r'')$ contradicting the fact that $f$ is one-to-on. $(2)$.  If $g_i(p,q)=g_j(p,q)$ for some pair $(p, q)$ then from Equation (\ref{eq:invarianceequality}) we get $g_i(p,q)=g_j(p,q)$ for all pairs $(p,q)$.
$(3)$. From  the  partition theorem of Erd\"os and Rado \cite{erdos-rado} (about partitions  of pairs of $\NN$ into an infinite numbers of colors) we may deduce that either $g_i$ is one-to-one or it does not depends on the first or it does not depends on the second coordinate or is constant. From this, the result follows. To give an hint of this deduction, we prove that $g_i$ is constant if there are $p<q<r$ such that $g_i(p, q)= g_i(q,r)$. Indeed,  let $p<q<r$ and $p'<q'<r'$ in $\NN$. Let $h$ be a local isomorphism of $(\mathbb N, \leq)$ sending $p,q, r$ onto $p',q', r'$.  If $g_i(p,q)=g_i(q,r)$, then applying Equation (\ref{eq:invarianceequality}), we get  $g_i(p',q')=g_i(q',r')$.  In particular, we get  $g_i(p,q)=g_i(q,\ell)$ for every $\ell>q$. Hence for every $p'<q'$, the value $g_i(p',q')$ does not depends upon $q'$.  With the fact that $g_i(p',q')=g_i(q',r')$ it follows that $g_i$ is constant. $(4)$. Set $h(f(p)):= f(p')$, and  $h(g_i(p,q)):= g_i(p',q')$.  Due to $(1)$, $(2)$ and $(3)$ and the invariance property, this  defines an isomorphism $h$ from $R_{\restriction \{f(p) \cup \mathcal F(p,q)\}}$ to $R_{\restriction \{f(p') \cup \mathcal F(p',q')\}}$.  Similarly, we may define an isomorphism from $R_{\restriction \{f(r) \cup \mathcal F(p,q)\}}$ to $R_{\restriction \{f(r') \cup \mathcal F(p',q')\}}$. The conclusion follows. \hfill $\Box$

\begin{definition}We say that $\mathfrak {L}$ has type $(I)$ if  the restrictions $R_{\restriction \{f(p) \cup \mathcal F(p,q)\}}$ and  $R_{\restriction \{f(r) \cup \mathcal F(p,q)\}}$ are isomorphic for some integers $p<q<r$.  Otherwise it has type $(II)$.
\end{definition}

\noindent \begin{claim}\label{keyclaim}Let $p<q<r$.
If $\mathfrak{L}$ has type $(I)$ then  the restrictions $R_{\restriction \{f(q) \cup \mathcal F(p,q)\}}$ and  $R_{\restriction \{f(r) \cup \mathcal F(p,q)\}}$ are not isomorphic, whereas if $\mathfrak {L}$ has type $(II)$
the restrictions $R_{\restriction \{f(p) \cup \mathcal F(p,q)\}}$ and  $R_{\restriction \{f(r) \cup \mathcal F(p,q)\}}$ are not isomorphic.
\end{claim}

\noindent{\bf Proof of Claim \ref{keyclaim}}.
Since $\mathfrak L$ is invariant, if  it has type $(I)$ then the restrictions $R_{\restriction \{f(p) \cup \mathcal F(p,q)\}}$ and  $R_{\restriction \{f(r) \cup \mathcal F(p,q)\}}$ are isomorphic. By construction,  the restrictions $R_{\restriction \{f(p) \cup \mathcal F(p,q)\}}$ and  $R_{\restriction \{f(q) \cup \mathcal F(p,q)\}}$ are not isomorphic, hence the restrictions $R_{\restriction \{f(q) \cup \mathcal F(p,q)\}}$ and  $R_{\restriction \{f(r) \cup \mathcal F(p,q)\}}$ are not isomorphic. By definition, if $\mathfrak L$ has type $(II)$ the restrictions $R_{\restriction \{f(p) \cup \mathcal F(p,q)\}}$ and  $R_{\restriction \{f(r) \cup \mathcal F(p,q)\}}$ are not isomorphic.
\hfill  $\Box$
\medskip

Next, we built our almost-chainable restriction. Let $F':=\{i: g_i \;  \text {is constant} \; \}$,  $F$ be the image of these $g_i$'s and $K':= \{0, \dots, k\} \setminus F'$. We may suppose that $F'$ is made of the last elements of $\{1, \dots k\}$, hence $K':= \{0,\dots k'\}$. Let $V':= (\mathbb N\times K') \cup F$. We define a map $\tilde \Phi$ from $V'$  to $V$ as follows.

We set $\tilde \Phi(x)=x$ for $x\in F$ and $\tilde \Phi(n,i):= g_i(2n, 2n+1)$ for $i=1,\dots, k'$. If $\mathfrak {L}$ has type $(I)$, we set  $\tilde \Phi(n,0):= f(2n+1)$  and if $\mathfrak {L}$ has type $(II)$, we set  $\tilde \Phi(n,0):= f(2n)$.   Let $range (\tilde \Phi)$ be the image of  $\tilde \Phi$.

\begin{claim}\label{lem:infiniteclasses}
The restriction of $R$ to $range(\tilde \Phi)$ is almost-multichainable and belongs to $\mathscr S_{\mu,k}$.
\end{claim}
\noindent{\bf Proof of claim \ref{lem:infiniteclasses}. }
The map $\tilde \Phi$ is one-to-one.  To see that, note first that the images of $F$ and $\NN\times K'$ are disjoint. Indeed, suppose otherwise that  $g_i(2n, 2n+1)\in F$ for some $i\in \{0, \dots, k'\}$ and $n\in \NN$.  Then, for some $j\in F'$, we have $g_i(2n, 2n+1)= g_j(2n, 2n+1)$.  If $i\not =0$, then according to $(2)$ of Claim \ref{basicclaim}, we have $g_i=g_j$ hence $g_i$ is constant, and thus $i\in F'$, contradicting our choice of $i$. If $i=0$, we have $f(m)= g_0(2n, 2n+1)= g_i(2n, 2n+1)$ with $m=2n+1$ if $\mathfrak{L}$ has type $(I)$ and $m=2n$ if $\mathfrak{L}$ has type $(II)$. But, by construction, neither $f(2n)$ nor $f(2n+1)$ belongs to $\{g_i(2n, 2n+1): i=1, \dots k'\}$.  To conclude, it suffices to show that $\tilde \Phi$ is one-to-one on $\NN\times K'$. Suppose that $\tilde \Phi(n,i)=\tilde \Phi(n',i')$. If $i=i'=0$ this yields $f(2n+1)= f(2n'+1)$ or $f(2n)= f(2n')$ depending on the type of $\mathcal L$. In both cases, we have $n=n'$ since $f$ is one-to-one. If  $i$ and $i'$ are distinct from $0$  we have $g_i(2n, 2n+1)= g_{i'}(2n', 2n'+1)$. If $i=i'$, then according to $(3)$ of Claim \ref{basicclaim},  
$g_i$ is constant contradicting the choice of $i$. If $i\not =i'$ then, since $g_i\not =g_{i'}$, we get $n\not =n'$ from $(2)$  of Claim \ref{basicclaim}. With $n''$ distinct from $n$ and $n'$,  we get, using equation \eqref{eq:invarianceequality}, $g_{i'}(2n',2n'+1)=g_{i'}(2n'',2n''+1)$. By $(3)$ of Claim \ref{basicclaim}, $g_{i'}$ is constant. Impossible. It remains to consider the case $i=0$, $i'\not =0$. It does not happen since by construction, and from $(1)$ of claim \ref{basicclaim}, $f(n)$ is not equal to an element of the form  $g_{i}(p,q)$.

Next, let  $S$ be the inverse image of $R$ by  $\tilde \Phi$. We need to prove that  for every local isomorphism  $h$ of $C:=(\mathbb N, \leq)$ the map $\hat h:=(h, 1_{K'})$
extended by the identity on $F$ is a local isomorphism of $S$ (the map
$(h, 1_{K'})$ is defined by $(f, 1_{K'})(x, y):= (f(x), y)$). This fact follows immediately from the invariance of $\mathfrak L$.  Now,  observe that
by construction, $(n,0) \not \simeq_{\leq k,S}(n',0)$ for $n<n'$. Indeed, set $A'_{n}:=(\{n\}\times K'\setminus\{0\}) \cup F$ and $n<n'$. We claim that $S_{\restriction \{(n,0)\} \cup A'_n}$ and $S_{\restriction \{(n',0)\} \cup A'_n}$ are not isomorphic. Indeed, let $A_n:=\tilde \Phi(A'_n)$.  The restriction  $S_{\restriction \{(n,0)\} \cup A'_n}$ is isomorphic to  $R_{\{\tilde \Phi(n, 0)\}\cup A_n}$ and the restriction $S_{\restriction \{(n',0)\} \cup A'_n}$ is isomorphic to  $R_{\{\tilde \Phi(n', 0)\}\cup A_n}$. We have $A_n:=\{ g_i(2n, 2n+1): i=1,\dots, k'\} \cup F$. If  $\mathfrak L$ has type $(I)$, $\tilde \Phi(n, 0)=f(2n+1)$,  $\tilde \Phi(n', 0)=f(2n'+1)$    and by Claim \ref{keyclaim}, $R_{\restriction \{f(2n+1)\} \cup A_n} \not \simeq R_{\restriction \{f(2n'+1)\} \cup A_n}$, proving that $(n,0) \not \simeq_{\leq k,S}(n',0)$. If  $\mathfrak {L}$ has type $(II)$, we have  $\tilde \Phi(n,0)= f(2n)$ and $\tilde \Phi(n', 0)=f(2n')$.  Again by Claim \ref{keyclaim}, $R_{\restriction \{f(2n)\} \cup A_n} \not \simeq R_{\restriction \{f(2n')\} \cup A_n}$, proving that $(n,0) \not \simeq_{\leq k,S}(n',0)$. Since all elements of the form $(n, 0)$ for $n\in \NN$ are not $\simeq_{\leq k, S}$-equivalent, it follows from Corollary \ref{lem:test} that  $S\in \mathscr S_{\mu,k}$, hence  $R_{\restriction range(\tilde \Phi)}\in \mathscr S_{\mu,k}$. This proves our claim. \hfill $\Box$

Since $S$ has type at most $(\omega, k+1)$, the conclusion of Theorem \ref{mainthm} follows.
\end{proof}

\subsection{Proof  of Theorem \ref{theo:base}}

 Let $R\in \mathscr O_{\mu}$. Let $m$ be the maximum of $\mu$. Let $k:=i(m)$. By Theorem \ref {thm:marytreshold}, $\simeq _{R}$ coincides with $\simeq _{\leq k,R}$. Hence, $R\in  \mathscr O_{\mu,k}$. Let $\mu'$ be the sequence obtained by putting $2$ in front of $\mu$. Then $R\in \mathscr S_{\mu', k}$.  Apply Theorem \ref{mainthm}.

\subsection{Proof  of Theorem \ref{theo:expo}}
Let  $R\in \mathscr O_{\mu, k}$.  We have to show that the profile of $R$ is bounded from below by a generalization of the Fibonacci sequence. According to Theorem \ref{mainthm}, $R$ embeds some ordered almost-multichainable relational structure  belonging to  $\mathscr S_{\mu',k}$. Thus it suffices to prove that the profile of $R'$ is bounded from below by a generalization of the Fibonacci sequence. Thus we may suppose as well that $R$ is multichainable and belongs to $\mathscr S_{\mu',k}$.

 Let $h$ be a positive  integer and let $(w_h(n))_{n\geq 0}$ be the sequence of integers defined by $w_h(n):=1$ for every $n$, $0\leq n<h$ and the recurrence relation $w_h(n)=w_h(n-1)+w_h(n-h)$ for $n\geq h$. For example, $w_2(n)=F(n+1)$ where $(F(m))_{m\in \mathbb N}$ is the Fibonacci sequence. Example for $h=1,\dots,   8$ are in OEIS (eg see A000045, A003269, A005710). 
The generating series  $w_h:=\sum_{n\geq 0}w_h(n)X^n$ satisfies

  \begin{equation}\label{eq:GF}
w_h= \frac{1}{1-X-X^h}.
 \end{equation}

 Asymptotically,  $w_h(n) \simeq d.c_{h}^{n}$ where $c_{h}$ is the largest positive solution of the characteristic equation $X^h-X^{h-1}-1=0$ and $d$ some positive constant.

 We are going to  find $s\leq k+1$ and  to enumerate some particular restrictions of $R$ and show that the number of such restrictions defined on $n$-element subsets verifies the recurrence relation above.  
 Doing so, it will follow that $\varphi_R(n)\geq w_s(n)$. 
 Asymptotically,  this  yields $\varphi_R(m)\geq d.c_{s}^m$. 
 This proves Theorem \ref{theo:expo}.


 Indeed,  according to Theorem \ref{mainthm} every member of $\mathscr O_{\mu, k}$ embeds some  almost-multichainable $R^{\prime}\in \mathscr O_{\mu, k}$. Applying this to $R$ and using the fact that  $c_{k+1}\leq c_{s}$ we get $\varphi_R(m)\geq d.c_{k+1}^m$ as claimed.

 According to Theorem  \ref{mainthm}, we may suppose that $R$ is of the form
$R:= ((\mathbb N\times K) \cup F, \leq,  (\rho_i)_{i\in I})$,  $\vert K\cup F\vert =k+1$, $K$  is a $k'+1$-element interval $[k'_K,k''_K]$ of $\mathbb Z$  containing $0$ such that all  elements $(n,0)$  for $n\in \mathbb N$ are inequivalent modulo $\simeq_R$. More precisely, for $n\in \NN$, set $A_{n}:= \{n\}\times  (K\setminus \{0\}) \cup F$; according to the proof of Claim \ref{lem:infiniteclasses},   we may suppose that for each $n\in \mathbb N$, the restrictions $R_{\restriction {\{(n, 0)\} \cup A_{n}}}$ 
and $R_{\restriction {\{(n+1, 0)\} \cup A_{n}}}$ are not isomorphic.  With no loss of generality we may suppose that the order induced by $\leq$ on $\{n\}\times K$ agree with the natural order; furthermore we may suppose that its restriction to $\mathbb N\times \{0\}$   also agree  with the natural order (if it disagree, we may associate a relational structure $R'$ to the natural order; it will have the same age as $R$,   hence the same profile). Let $r:= k-k'= \vert F\vert $ and $h:=k'+1$.
The chain $C:= (\mathbb N\times K, \leq)$ decomposes into three intervals $A', A $ and $A''$ as follows. Let $K_1$ be the largest interval of $[k'_K, k''_K]$ containing $0$ such that $(0,Max (K_1))\leq (1, Min (K_1))$.  Set $k'_1:= Min (K_1)$, $k''_1:= Max (K_1)$, $A':= \mathbb N\times [k'_K,k'_1-1]$, $A:= \mathbb N\times [k'_1,k''_1]$,  $A'':= \mathbb N\times [k''_1+1, k''_K]$, $s':=\vert [k'_K,k'_1-1]\vert$, $s:= \vert [k'_1,k''_1]\vert=\vert K_1\vert$ and $s'':= \vert [k''_1+1, k''_K]\vert$. Due to the invariance property of $C$, $A'$ is an initial interval,  $A''$ is  a final interval and $C$ decomposes into the lexicographical sum $A'+A+A''$; furthermore $A$ decomposes into the lexicographical sum $\sum_{n\in \mathbb N} \{n\}\times K_1$. The invariance property also implies that each element of $F$ is either less than every element of $A'$ or greater than every element of $A''$.
Set $F'$ the set of elements of $F$ which are less than every element of $A'$ and  $F''$ those which are greater than the elements of $A''$, $r':= \vert F'\vert$ and $r'':= \vert F''\vert$. We have $F=F'\cup F''$ hence $r:=r'+r''$. The chain $((\mathbb N\times K) \cup F, \leq )$ decomposes into the sum 
$(F'+A')+A+(A''+F'')$. 

To each integer $n$, we will enumerate all non isomorphic restrictions of $R$ defined on $n$-element subsets of the form $E:= (I\times \{0\})\cup (J\times K_1)$ where $I\cup J$ is a finite subset of $\NN$ such that $\vert I\vert+s\vert J\vert=n$.

We recall that the chain $C_{\restriction E}$ is such that $(m,i) < (m',i')$ if $(m<m')$ or $(m=m' \text{ and } i < i')$.

Set, for each $n\in\NN$, $w_s(n)$ the number of non isomorphic restrictions of $R$ on $n$-element subsets of the form $E$. If $n<s$, then $J=\varnothing$ and the set is of the form $I\times\{0\}$. Hence $w_s(n):=1$. If $n\geq s$, then we consider the $n$-elements ordered w.r.t $\leq$. By invariance, the restrictions of $R$ to the sets $\{m\}\times K_1$ with $m\in\NN$ are all isomorphic. Then, either the chain ends by the set $\{m\}\times K_1$ with $m\in J$, and, in this case, the number of non isomorphic substructures with $n$-elements of the form $E$ is $w_s(n-s)$ or not. And in this later case, the chain ends by an elements $(m,0)$ with $m\in I$ and the number of such substructures of order $n$ is $w_s(n-1)$. We then get
$$\left\{\begin{array}{l}
w_s(n)=1 \;\text{ for }0\leq n<s.\\
w_s(n)=w_s(n-s)+w_s(n-1)\;\text{ for }n\geq s.\end{array}
\right.$$ 

With this, the proof of Theorem \ref{theo:expo} is complete. 

\section{The case of graphs}\label{section:graphs-ordered}

The graphs we consider in this section are undirected without loops.
We recall some properties of autonomous subsets of graphs (see \cite{ehrenfeucht, schmerl-trotter}). Let $G$ be a graph, we  recall that a subset $A\subseteq V(G)$ is \emph{autonomous} if for every $x, x'\in A$, $y\notin A$, the following condition holds:$$\{x,y\}\in E(G) \text{ if and only if } \{x',y\}\in E(G).$$
The empty set, the one-element subsets of $V(G)$ and the hole set $V(G)$ are autonomous. If $G$ has no other autonomous subset,  $G$ is said \emph{indecomposable}; if $G$ has more than two elements it contains a $P_4$ (a path on four vertices); in this case $G$ is \emph{prime}, this is a well known result of Sumner \cite{sumner}.

We also recall that if a graph $G$ is a lexicographical sum $\underset{i\in H}\sum G_i$ of a family of graphs $G_i$, indexed by a graph $H$, then, provided that they are pairwise disjoint, the $V(G_i)$ form a partition of $V(G)$ into autonomous sets. Conversely, if the set $V(G)$ is partitioned into autonomous sets, then $G$ is the lexicographical sum of the graphs induced on the blocks of the partition.

With these definitions, we have immediately:

\begin{lemma} \label{lem-one-equiv}Let $G:= (V(G), E(G))$ be a undirected graph. Then
the equivalence relations $\simeq _{1, G}$ and $\simeq _{G}$ coincide. In this case,
     equivalence classes are maximal autonomous subsets which are cliques or independent sets.
\end{lemma}

We will obtain the ten graphs mentioned  in  Proposition \ref{prop: graphs}  via Ramsey's theorem under the form of Theorem \ref{thm:ramsey-invariant}. These ten graphs are almost-multichains on $F\cup(L\times K)$ such that $L:=\NN$, $\vert K\vert=2$ and $F=\emptyset$. 
Before describing these graphs, let's give the way to construct them.


\smallskip

 Let $G$ be a graph which does not decompose into a finite lexicographical sum $\underset{i\in H}\sum G_i$ of a family of undirected graphs $G_i$,  each one being a clique or an independent set, indexed by a finite undirected graph $H$. Then the equivalence relation $\simeq_G$  has  infinitely many equivalence classes. Thus,  there is a one-to-one map $f:\NN\longrightarrow V(G)$ such that the images of two distinct elements are pairwise inequivalent. Due to Lemma \ref {lem-one-equiv}, we may  define a map $g:[\NN]^2\longrightarrow V(G)$ such that $g(n,m)$ witnesses the fact that $f(n)$ and $f(m)$ are inequivalent for every $n<m$. That is  $$\{f(n),g(n,m)\}\in E(G)\Leftrightarrow\{f(m),g(n,m)\}\notin E(G).$$

Let $\omega:=(\NN,\leq)$, $\Phi:=\{f,g\}$ and $\mathfrak{L}:=\left\langle \omega,G,\Phi \right\rangle $. Ramsey's theorem allows us to find an infinite subset $X$ of $\NN$ such that $\mathfrak L\restriction_X$ is invariant.    

By relabeling $X$ with integers, we may suppose $X=\NN$, thus $\mathfrak{L}$ is invariant.

\begin{claim}\label{graphclaim}
\begin{enumerate}
\item $\{f(n),f(m)\}\in E(G)\Leftrightarrow\{f(n'),f(m')\}\in E(G),~\forall n<m,~n'<m'$.
\item $\{g(n,m),g(n',m')\}\in E(G)\Leftrightarrow\{g(k,l),g(k',l')\}\in E(G),~\forall n<m<n'<m',~k<l<k'<l'$.
\item $\{f(n),g(n,m)\}\in E(G)\Leftrightarrow\{f(k),g(k,l)\}\in E(G),~\forall n<m,~k<l$.
\item $\{f(m),g(n,m)\}\in E(G)\Leftrightarrow\{f(l),g(k,l)\}\in E(G),~\forall n<m,~k<l$.
\item $\{f(k),g(n,m)\}\in E(G)\Leftrightarrow\{f(l),g(p,q)\}\in E(G),~\forall n<m<k,~p<q<l$.
\item $\{f(n),g(m,k)\}\in E(G)\Leftrightarrow\{f(p),g(q,l)\}\in E(G),~\forall n<m<k,~p<q<l$.
\item If  $\{f(n),g(n,m)\}$ and  $\{f(k),g(n,m)\}$ are both edges or non edges for some integers $n<m<k$ then $\{f(n'),g(n',m')\}$ and  $\{f(k'),g(n',m')\}$ are both edges or non edges for every $n'<m'<k'$.
\item $g(n,m)$ and $f(k)$ are different for every integers $n<m$ and $k$ .
\item $g(n,m)\neq g(n',m')$ for every $n<m<n'<m'$.
\end{enumerate}
\end{claim}

\textbf{Proof of Claim \ref{graphclaim}.}
The seven  first items follow from invariance. To prove $(8)$, suppose that there are integers $n<m<k$ such that $g(n,m)=f(k)$. By construction of the functions $f$ and $g$, we have $\{f(n),g(n,m)\}\in E(G)\Leftrightarrow\{f(m),g(n,m)\}\notin E(G)$. According to $(1)$, $\{f(n),f(k)\}$ and $\{f(m),f(k)\}$ are both edges or non edges, replacing $f(k)$ by $g(n,m)$ we get a contradiction with the above equivalence. Item $(9)$. Suppose that there are integers $n<m<n'<m'$ such that $g(n,m)=g(n',m')$. According to $(5)$, 
$\{f(n'),g(n,m)\}$ and  $\{f(m'),g(n,m)\}$ are both edges or non edges. Replacing $g(n,m)$ by $g(n',m')$, we get $\{f(n'),g(n',m')\}$ and $\{f(m'),g(n',m')\}$ are both edges or non edges. This contradicts the above equivalence.

\hfill       $\Box$

 Let $H:\NN\times\{0,1\}\longrightarrow V(G)$. Set $H(n,1):=g(2n,2n+1)$. From $(7)$ 
 of Claim \ref{graphclaim} we have two cases. If $\{f(n),g(n,m)\}$ and  $\{f(k),g(n,m)\}$ are both edges or non edges for some integers $n<m<k$, then set $H(n,0):=f(2n+1)$. Otherwise  set $H(n,0):=f(2n)$.

Let $G'$ be the undirected graph with vertex set $V(G'):=\NN\times\{0,1\}$ such that $$\{x,y\}\in E(G')\Leftrightarrow\{H(x),H(y)\}\in E(G)$$
for every pair $x,y$ of vertices of $\NN\times\{0,1\}$. By construction of $G'$ we have:
\begin{equation}\label{eq1}
\quad\{(n,0),(n,1)\}\in E(G')\Leftrightarrow\{(m,0),(n,1)\}\notin E(G')
\end{equation}
for every integers $n<m$, then $(n,0)$ and $(m,0)$ are inequivalent for every $n<m$ and hence $G'$ has infinitely many equivalence classes.
 By invariance of $\mathfrak L$ the graph $G'$ is an almost-multichain and then $G\restriction_{range(H)}$ is almost-multichainable.\\ 


 Now, to construct $G'$ it suffices to decide on relations between the four vertices $(0,0)$, $(0,1)$, $(1,0)$ and $(1,1)$, the remaining relations will be deduced by taking local isomorphisms of $C:=(\NN,\leq)$. For example if $\{(0,0),(0,1)\}\in E(G')$, then $\{(0,1),(1,0)\}\notin E(G')$.  We have then to decide on the relations $\{(0,0),(1,0)\}$, $\{(0,0),(1,1)\}$  and $\{(0,1),(1,1)\}$, this leads to eight ($2^3$) cases but two of these graphs are isomorphic, thus  we have only seven cases. If  $\{(0,0),(0,1)\}\notin E(G')$, then $\{(0,1),(1,0)\}\in E(G')$;  we have then eight other cases but five have already been obtained, hence we get only three new graphs. The total is ten. \\
Note that,  if $\{(0,0),(1,0)\}\in E(G')$ then $\{(n,0),(m,0)\}\in E(G')$ for every $n<m$ and thus $\{(n,0), ~n\in\NN\}$ is a clique and if $\{(0,0),(1,0)\}\notin E(G')$ then $\{(n,0),(m,0)\}\notin E(G')$ for every $n<m$ and thus $\{(n,0), ~n\in\NN\}$ is an independent set.
\smallskip

\subsection{Description of the ten graphs.}

Denote by $G_i, 1\leq i\leq 10$, the ten graphs defined on the same set of vertices $V:=\mathbb N\times\{0,1\}$. Set  $A:=\mathbb N\times\{0\}$ and $B:=\mathbb N\times\{1\}$. For $i:=1,2,3$, the subsets $A$ and $B$ are independent sets and a pair  $\{(n,0),(m,1)\}$, $n,m\in\mathbb N$ is an edge in $G_1$ if $n=m$, an edge in $G_2$ if $n\leq m$ and an edge in $G_3$ if $n\neq m$. Thus,  $G_1$ is a direct sum of infinitely many copies of $K_2$ (the complete graph with two vertices) and $G_2$ is the half complete bipartite graph, a critical graph of Schmerl and Trotter \cite{schmerl-trotter}. For $4\leq i\leq 7$, one of the subsets  $A$ or $B$ is a clique and the other is an independent set. Thus, the set $E(G_4)$ of edges of  $G_4$ is $E(G_1)\cup\{\{x,y\}, x\neq y\in A\}$, 
 $E(G_5):=E(G_2)\cup\{\{x,y\}, x\neq y\in A\}$, 
$E(G_6):=E(G_2)\cup\{\{x,y\}, x\neq y\in B\}$ 
and $E(G_7):=E(G_3)\cup\{\{x,y\}, x\neq y\in A\}$. 
Thus, in the graphs $G_4$, $G_5$, and $G_7$ the set $A$ is a clique, $B$ is an independent set and these graphs coincide with $G_1$, $G_2$, and $G_3$ respectively on pairs $\{x,y\},~x\in A, y\in B$. For $G_6$, the set $B$ is a clique, $A$ an independent set and it coincides with $G_2$ on pairs $\{x,y\},~x\in A, y\in B$.
The graphs  $G_8$, $G_9$ and $G_{10}$ are such that the subsets $A$ and $B$ are both cliques, they coincide with $G_1$, $G_2$, and $G_3$ respectively on pairs $\{x,y\},~x\in A, y\in B$.
The graph $G_8$ is the dual of $G_3$, the graph $G_7$ is the dual of $G_4$,  the graph $G_{10}$ is the dual of $G_1$, the graphs $G_5$ and $G_6$ are equimorphic, but non isomorphic, with their dual respectively and have the same ages, each of the graphs $G_2$ and $G_9$ embeds the dual of the other.

\subsection{Profiles of the ten graphs.} We call \emph{profile} of a graph $G$ and denote by $\varphi_G$ the profile of the age of   $G$. Clearly, a graph and its dual have the same profile and if a graph $G$ embeds some graph $H$ then $\varphi_H\leq \varphi_G$. Then, according to the description above, the graphs $G_1$, $G_2$, $G_3$, $G_4$ and $G_5$ have the same profiles as $G_{10}$, $G_9$, $G_8$, $G_7$ and $G_6$ respectively. From the analysis below, $\varphi_{G_1}$, $\varphi_{G_3}$ and $\varphi_{G_4}$ have polynomial growth and $\varphi_{G_2}$, $\varphi_{G_5}$ have exponential growth.

 \textbf{Profile of $G_1$:} We have $\varphi_{G_1}(n)=\lfloor\frac{n}{2}\rfloor +1,\forall n\in\mathbb N$. Indeed, every subgraph of $G_1$ defined on $n$ vertices is isomorphic, for some integers $p, q$, to $pK_2\oplus q$,  which is the direct sum of $p$ copies of $K_2$, the clique with two vertices, and an independent set of size $q$, such that $p\leq \frac{n}{2}$ and $q=n-2p$. Then we may represent every subgraph of size $n$ by a couple $(n,p)$ with $p\leq \frac{n}{2}$ and inversely. Two subgraphs represented by $(n,p)$ and $(n',p')$ are isomorphic iff $n=n'$ and $p=p'$. Its generating series is the rational function $$F_{G_1}(x)=\dfrac{1}{(1-x)(1-x^2)}.$$

\textbf{Profile of $G_2$:} The first values are $1,1,2,3,6,10,20,36,72,136$. It satisfies the following recurrence: $\varphi_{G_2}(n)=\varphi_{G_2}(n-1)+2^{n-3}$ for $n\geq 3$ and $n$ is odd; $\varphi_{G_2}(n)=\varphi_{G_2}(n-1)+2^{n-3}+2^{\frac{n-4}{2}}$ for $n\geq 4$ and $n$ even. Thus $\varphi_{G_2}(n)\geq 2^{n-3}$ for $n\geq 3$. Its generating series is the rational function $$F_{G_2}(x)=\dfrac{1-x-2x^2+x^3}{(1-2x)(1-2x^2)}.$$ 

\textbf{Profile of $G_3$:} The first values are $1,1,2,3,6,6,10,10$, and \\$\varphi_{G_3}(n)=\underset{k=0}{\overset{\lfloor n/2\rfloor}{\sum}}(k+1)$ for every $n\neq 2$. Its generating series is the rational function $$F_{G_3}(x)=\dfrac{1-x^2+x^3+2x^4-2x^5-x^6+x^7}{(1-x)(1-x^2)^2}.$$ 

\textbf{Profile of $G_4$:} The first values are $1,1,2,4,7,10,14,18,23,28$. It satisfies the following relation: $\varphi_{G_4}(n)=\dfrac{1}{4}(n-1)(n+5)$ for $n\geq 3$ and $n$ is odd; $\varphi_{G_4}(n)=\frac{1}{4}n(n+4)-1$ for $n\geq 4$ and $n$ even. Its generating series is the rational function $$F_{G_4}(x)=\dfrac{1-x+2x^3-x^5}{(1-x)^3(1+x)}.$$ 

\textbf{Profile of $G_5$:} $\varphi_{G_5}(n)=2^{n-1}$ for $n\geq 1$. The first values are $1,1, 2, 4,8, 16, 32$. Its generating series is the rational function $$F_{G_3}(x)=\dfrac{1-x-2x^2}{(1-2x)}.$$
\bigskip

\end{document}